\documentclass{elsart}
\usepackage{stmaryrd}%bigsqcup
\usepackage{amsfonts}
\usepackage{amsmath}
\usepackage{mathrsfs}
\usepackage{xcolor}

\def\dda{\mathord{\mbox{\makebox[0pt][l]{\raisebox{-.4ex}
{$\downarrow$}}$\downarrow$}}}

\newtheorem{tm}{Theorem}[section]
\newtheorem{pn}[tm]{Proposition}
\newtheorem{lm}[tm]{Lemma}

\theoremstyle{definition}
\newtheorem{dn}[tm]{Definition}
\newtheorem{rk}[tm]{Remark}
\newtheorem{ex}[tm]{Example}

\journal{Elsevier}
\usepackage{bm}

\begin{document}
	
\begin{frontmatter}
\title{Representations of
domains via
closure spaces in the quantale-valued setting}
\author{Guojun Wu$^{1,2}$, Wei Yao$^{1,2}$, Qingguo Li$^{3}$}
\address{$^{1}$School of Mathematics and Statistics, Nanjing University of Information Science and Technology, Nanjing, 210044, China\\
$^{2}$Applied Mathematics Center of Jiangsu Province, Nanjing University of Information Science and Technology, Nanjing, 210044, China\\
$^{3}$School of Mathematics, Hunan University, Changsha, 410082, China
}
\date{}

\begin{abstract}
With a commutative unital quantale
$L$ as the truth value table, this study focuses on the representations of
$L$-domains by means of  $L$-closure spaces. First, the notions of interpolative generalized
$L$-closure spaces and directed closed sets are introduced. It is proved that in an interpolative generalized
$L$-closure space (resp.,
$L$-closure space),  the collection of directed closed sets with respect to the inclusion
$L$-order forms a continuous
$L$-dcpo (resp., an algebraic
$L$-dcpo). Conversely, it is shown that every continuous
$L$-dcpo (resp., algebraic
$L$-dcpo) can be reconstructed by an interpolative generalized
$L$-closure space (resp.,
$L$-closure space). Second, when
$L$ is integral, the notion of dense subspaces of generalized
$L$-closure spaces is introduced. By means of dense subspaces, an alternative representation for algebraic
$L$-dcpos is given. Moreover, the concept of
$L$-approximable relations between interpolative generalized
$L$-closure spaces is introduced. Consequently, a categorical equivalence between the category of interpolative generalized $L$-closure spaces (resp., $L$-closure spaces) with $L$-approximable
relations and that of continuous  $L$-dcpos (resp., algebraic  $L$-dcpos) with Scott continuous mappings is established.

%With a commutative unital quantale $L$ as the truth value table, we investigate the
%representations for $L$-domains using  $L$-closure spaces. Firstly, we introduce the notions
%of interpolative generalized $L$-closure spaces and directed closure sets. We prove that the collection of directed closure sets in
%an interpolative generalized $L$-closure space (resp., $L$-closure space)
%with respect to inclusion $L$-order is a continuous $L$-dcpo (resp., an algebraic $L$-dcpo). Conversely, we show that every
%continuous $L$-dcpo (resp.,  algebraic $L$-dcpo) can be
%captured by an interpolative generalized $L$-closure space (resp., $L$-closure space). Secondly, when $L$ is integral,
%we introduce the notion of dense subspaces of generalized $L$-closure spaces. By means of dense subspaces,
%we give another representation for algebraic $L$-dcpos. Moreover, we introduce the notion of $L$-approximable
%relations between interpolative generalized $L$-closure spaces.  Consequently, we establish a categorical equivalence
%of interpolative generalized $L$-closure spaces (resp., $L$-closure spaces) with $L$-approximable
%relations and that of continuous  $L$-dcpos (resp., algebraic  $L$-dcpos) with Scott continuous mappings.
%
%

\end{abstract}

\end{frontmatter}

\section{Introduction}

Domain theory \cite{Gierz,Jean-2013} is a multidisciplinary field of mathematics that studies a specific class of ordered structures known as domain structures. These structures have extensive applications in theoretical computer science, serving as denotational models for functional languages.
Research in domain theory dates back to  to the pioneering work of  Scott \cite{Scott-72} in the  1970s. A fundamental feature of domain theory is the interaction among various mathematical structures, such as orders, algebras, topologies, and logics.
The representation theory of various domains is a crucial topic that highlights the connections between ordered structures and other mathematical
structures. The representation of domains refers to finding a
concrete and accessible alternative for abstract structures of domains \cite{Wang-Li-studia}.
The family of sets with respect to  the set-theoretic inclusion
is the simplest ordered structure, inspiring scholars to use appropriate set-families, induced by specific mathematical structures, to represent various domains.  These mathematical structures include topological spaces \cite{Wang-Li-studia,XU-mao-form,XU-Yang}, information systems \cite{He-Xu,Spreen-Xu,Wang-Li-semi} and non-classical logics \cite{Wang-log-2024,Wang-log1}.

A closure space consists of a set
and a closure operator on it.  There are close links between closure spaces and ordered structures.
Many studies have shown that closure spaces are closely related to  domain theory and  play a key role in various completions of
ordered structures (see \cite{Day,Zhao, ZX-Zhang, Zhang-Shi-Li-2021}).
Recently, an increasing number of studies
have shown that closure spaces are useful tools for representing  various domains.
 The idea of representing ordered  structures by  closure spaces with some additional conditions can be traced back to Birkhoff's
representation theorem for finite distributive lattices \cite{Brikhoff} and Stone's duality for
Boolean algebras \cite{Stone}.  In \cite{Erne-stone},  Ern\'{e}  demonstrated how various Stone-type
dualities arise from that general construction scheme. In particular,
 Ern\'{e} established  a Stone-type duality between sober algebraic closure spaces and  algebraic lattices. This work
 provided a representation of algebraic lattices via algebraic closure spaces.
Convex spaces are special algebraic closure spaces with the empty set as a closed set \cite{Van}.
In  \cite{Yao-Zhou}, Yao and Zhou established a categorical isomorphism between   join-semilattices and  sober convex spaces.
 In \cite{shen-zhao-shi}, Shen et al. established   a Stone-type duality between  sober convex spaces and algebraic lattices; moreover, they obtained
 additional results revealing the connection between closure
spaces and domain theory.

In order to use closure spaces to represent algebraic domains, Guo and Li generalized the notion of algebraic closure spaces to F-augmented closure spaces by adding a structure to a closure space, successfully providing representations for algebraic domains \cite{Guo-Li}.
Following Guo and Li's work, scholars simplified the structure of F-augmented closure spaces from different perspectives and used the resulting simplified closure spaces to represent algebraic domains \cite{Wu-Guo-Li, Wu-Xu-FIE}.
To achieve a transition from representing algebraic domains to representing continuous ones, Li et al. generalized the notion of  algebraic closure spaces to  continuous closure spaces and extended the representation theory from algebraic lattices to continuous lattices \cite{Li-Wang-Yao}. Recently, Wang and Li used interpolative generalized closure spaces to represent continuous domains, successfully establishing a categorical equivalence between interpolative generalized closure spaces with approximable mappings and domains with Scott continuous functions \cite{Wang-Li-studia}.

Quantitative domains are extensions of classical domains achieved by generalizing ordered relations to more general structures,
such as enriched categories, generalized metrics and lattice-valued  orders. This paper focuses on lattice-valued  orders.
The lattice-valued  order approach to quantitative domain theory , known as lattice-valued domain theory,  was initiated by  Fan
and Zhang \cite{FanTCS,FanZhang}, and further developed primarily  by Zhang et al. \cite{Yu-Zhang,zhang-gao}, Yao  \cite{PartI,YaoTFS}, Li et al \cite{LiQG1,Su-Li-Intel},
Zhao et al. \cite{Zhao2,Zhao3,Wang-Zhao}, Yue et al. \cite{Liu-Meng-Yue}, Guti\'{e}rrez Garc\'{i}a et al. \cite{Garcia}, among others.
Several well-known results in classical domain theory, which reflect the intersection of ordered structures and topological structures, have been extended to the lattice-valued setting. For instance, Yao \cite{YaoTFS} generalized Scott's result \cite{Scott-72} to the frame-valued setting,  establishing a categorical isomorphism between continuous $L$-lattices and injective $T_0$ $L$-topological spaces. Additionally,  Yao and Yue \cite{Yao-monad} extended Day's result \cite{Day} to the frame-valued setting,  providing  an algebraic representation of continuous $L$-lattices via
the open filter monad over the category of $T_0$ $L$-topological spaces.

Numerous studies have shown that some connections between closure spaces and ordered structures can also be technically extended to the lattice-valued setting.
 Inspired by Yao's work \cite{YaoTFS}, Xia \cite{Xia-FFS} extended Jankowski's work \cite{Jankowski} and established a categorical isomorphism between the category of balanced $L$-S${_0}$-convex spaces and that of  $L$-frames.
Su and Li generalized the Stone-type duality between algebraic lattices and algebraic closure spaces to the lattice-valued setting \cite{Su-Li-Intel}.
Yao and Zhou  established a  Stone-type duality between algebraic $L$-lattices and  sober $L$-convex spaces using a new approach \cite{Yao-Zhou-2021}.
Recently, Zhang \cite{ZX-Zhang,Zhang-Shi-Li-2021} discussed Z-completions of $\mathcal{Q}$-ordered sets
via the $\mathcal{S}$-completion 
of $\mathcal{Q}$-closure spaces 
within a more general framework, including a lattice-valued version of D-completion for $\mathcal{Q}$-cotopological as a special case.
As mentioned above, some connections between closure spaces and ordered structures can be technically extended to the lattice-valued setting.
It is thus natural to wonder  whether one can extend the categorial equivalence between interpolative closure spaces
and continuous domains  to the  lattice-valued  setting. In the present paper, we will examine  representations
for continuous $L$-dcpos, as well as algebraic $L$-dcpos, via generalized $L$-closure spaces.

 In this paper, we adopt the commutative unital quantale $L$  as the truth value table. The content of this paper is arranged as follows. In Section 2,
we recall basic concepts and results that will be frequently used throughout the
paper.
In Section 3, we  focus on the correspondence between interpolative generalized $L$-closure spaces and continuous  $L$-dcpos.
In Sections 4,  we focus on the correspondence between  $L$-closure spaces and algebraic $L$-dcpos.
In  Section 6, we introduce the
notion of  $L$-approximable relations between interpolative generalized closure
spaces, which leads to a categorically equivalent to that of continuous $L$-dcpos, as well as  algebraic $L$-dcpos, with Scott
continuous mappings as morphisms.

\section{Preliminaries}
We refer to \cite{residuated,Hajek,Quantale} for contents on quantales. We refer to \cite{MVTop} for notions of fuzzy sets,
and to \cite{YaoTFS,PartI} for contents on fuzzy posets. We refer to \cite{Well} for notions of  category.
\subsection{Commutative unital quantale, the truth value table}
%Let $L$ be a complete lattice and let $*$ be a semigroup operation
%on $L$. The pair $(L,*)$ is called a {\it quantale}
%if the operation $*$ is distributive over joins, that is,
%$$a*(\bigvee S)=\bigvee (a*S)\ \mbox{and}\ (\bigvee
%S)*a=\bigvee(S*a)$$ for all $a\in L$ and $S\subseteq L$, where $a*S=\{a*s|\ s\in S\}$ and $S*a=\{s*a|\ s\in S\}$. A
%quantale $(L,*)$ is called {\it commutative} (resp., {\it unital})
%if the operation $*$ is commutative (resp., has a unit element
%$k$).
%
%For a commutative quantale $(L,*)$, the operation
%$*$ has a right adjoint $\rightarrow:L\times L\longrightarrow L$
%given by $a\rightarrow b=\bigvee\{c\in L|\ a*c\leq b\}\ (\forall
%a,b\in L),$ satisfying that $$a*b\leq c\ \Leftrightarrow\ a\leq b\rightarrow c\ (\forall
%a,b,c\in L).$$

Let $L$ be a complete lattice with  a bottom element $0$ and a top element 1 and let $\otimes$ be a binary operation
on $L$.  The pair $(L,\otimes)$ is called a {\it quantale}
if the operation $\otimes$ is distributive over joins, that is,
$$a\otimes(\bigvee S)=\bigvee_{s\in S} (a\otimes s)\ \mbox{and}\ (\bigvee S)\otimes a=\bigvee_{s\in S} (s\otimes a).$$
for all $a\in L$ and $S\subseteq L$. A quantale  $(L,\otimes)$ is said to be {\em commutative} (resp., {\em unital}) if the operation is commutative (resp., has a unit element
$u$). A quantale  $(L,\otimes)$ is said to be  {\em integral} if it  is unital and the unit $u$ is the top element  of $L$, i.e., $u=1$.
For a  commutative  quantale $(L,\otimes)$, the operation
$\otimes$ gives rise to a right adjoint $\rightarrow:L\times L\longrightarrow L$
via the adjoint property  $$a\otimes b\leq c\ \Longleftrightarrow\ a\leq b\rightarrow c\ (\forall
a,b,c\in L).$$
The binary operation $\rightarrow$ is called the {\em implication}, or the {\em residuation}, corresponding to $\otimes $.

A quantale is called a
{\em frame}, or {\em complete Heyting algebra} if $\otimes=\wedge$.

%In this paper, unless otherwise stated, $L$ always
%denotes a complete residuated lattice with the operation $\otimes$.

Some basic properties of the binary operations $\otimes$ and $\rightarrow$ are listed below.
\begin{lm}\label{lm-resi-lat}{\rm (\cite{residuated,Hajek})} Suppose that $(L,\otimes, u)$ is a commutative unital quantale. Then for all $a,b,c\in L$ and  $\{a_i|\ i\in I\},\ \{b_j|\ j\in J\}\subseteq L$,

{\rm (Q1)} $u\leq a\rightarrow b\Longleftrightarrow a\leq
b$;

{\rm (Q2)} $ 0\rightarrow a=1$;

{\rm (Q3)} $u\rightarrow a=a$;

{\rm (Q4)} $a\otimes(a\rightarrow b)\leq b$;

{\rm (Q5)} $a\rightarrow (b\rightarrow c) = (a\otimes b)\rightarrow
c$.

{\rm (Q6)} $(\bigvee_{i\in I}a_i)\rightarrow b=\bigwedge_{i\in I}(a_i\rightarrow
b)$;

{\rm (Q7)} $a\rightarrow(\bigwedge_{j\in J}b_j)=\bigwedge_{j\in J}(a\rightarrow
b_j)$.
%
%{\rm (Q8)} $(c\rightarrow a)\rightarrow (c\rightarrow b)\geq
%a\rightarrow b$;
%
%{\rm (Q9)} $(a\rightarrow c)\rightarrow (b\rightarrow c)\geq
%b\rightarrow a$;

\end{lm}
{\bf Standing Assumption.} In this paper,   unless otherwise specified, the truth value table $L$ is always  assumed to be  a commutative  unital quantale, with unit denoted
by $u$.

\subsection{fuzzy sets and fuzzy posets}
Every mapping $A:X\longrightarrow L$ is called an $L$-{\it subset} of $X$ and we  use $L^X$ to denote the collection
of  $L$-subsets of $X$. Customarily, the crisp order $\leq$ on $L^X$ is  defined  pointwisely;
that is $A\leq B\Leftrightarrow A(x)\leq B(x)\ (\forall x\in X)$. An $L$-subset $A$ is said to be {\it nonempty} if $\bigvee_{x\in X}A(x)\geq u$.
Let $Y\subseteq X$ and $A\in L^X$, define $A|_{Y}\in L^Y$ by $A|_{Y}(y)=A(y)$ $(\forall y\in Y)$.
For an element $a\in L$, the notation $a_X$ denotes the constant $L$-subset of $X$ with  the value $a$, i.e., $a_X(x)=a\ (\forall x\in X)$.
For all $a\in L$ and $A\in  L^{X}$,  write $a\otimes A$ , for the $L$-subset given by $(a\otimes A)(x) = a\otimes A(x)$.

A {\em (binary) $L$-relation} $R$ from $X$ to $Y$ is a mapping
 $$R:X\times Y\longrightarrow L.$$
The {\em composition} $Q\circ R\in L^{X\times Z}$ of two $L$-relations $R\in L^{X\times Y}$ and $Q\in L^{Y\times Z}$ is  defined by
$$Q\circ R(x,z)=\bigvee_{y\in Y}R(x,y)\otimes Q(y,z)\ (\forall x\in X, z\in Z).$$
For each $a\in X$, define characteristic function  $u_a:X\longrightarrow L$    by

\centerline{
$
\ \ u_a(x)=\left\{\begin{array}{ll}u,& x=a;\\
0,&  x\neq a,
\end{array}\right.
$
}
 It is worth noting that in this paper, we also use the symbol $u_a$  to denote the characteristic function from $X^{\prime}$ to $L$. The only difference is the domains of the mappings.  To prevent confusion, readers can identify the domain of a given characteristic function based on the context.

Let $f:X\longrightarrow Y$ be a mapping between two sets. The {\it Zadeh extensions}
$f^\rightarrow:L^X\longrightarrow L^Y$ and $f^\leftarrow:L^Y\longrightarrow L^X$ are respectively given by
$$f^\rightarrow(A)(y)=\bigvee\limits_{f(x)=y}A(x)\ (\forall A\in L^X),\   \ \qquad f^\leftarrow(B)=B\circ f\ (\forall B\in L^Y).$$

\begin{dn}{\em (\cite{PartI})}
A mapping $e:P\times P\longrightarrow L$ is called
 an {\it $L$-order} if

{\rm (1)} $\forall x\in P$, $e(x,x)\geq u$;

{\rm (2)} $\forall x,y,z\in P$, $e(x,y)\otimes e(y,z)\leq e(x,z)$;

{\rm (3)} $\forall x,y\in P$, if $e(x,y)\wedge e(y,x)\geq u$,
then $x=y$.

The pair $(P,e)$ is called an $L$-{\em ordered set}. It is customary to write $P$ for the pair $(P, e)$.
\end{dn}
To avoid confusion, we sometimes use $e_P$ to denote the $L$-ordered set on $P$.
A mapping  $f: P\longrightarrow Q$  between two $L$-ordered sets is said to be {\it $L$-order-preserving}
 if for all $x, y\in P$, $e_{P}(x, y)\leq e_{Q}(f(x), f(y))$; $f$  is said to be  {\it $L$-order-isomorphic}
 if $f$ is a bijection and $ e_{P}(x, y)= e_{Q}(f(x), f(y))$ for all $x, y\in P$, denoted by $P\cong Q$.
\begin{ex}\label{ex-L-ord}{\em (\cite{YaoTFS})}
{\rm(1)}For a commutative unital quantale $L$,  define $e_L:L\times L\longrightarrow L$ by $e_L(x, y)= x\to y$ $(\forall x, y\in L)$.   Then $e_L$ is an $L$-order on $L$.

{\rm(2)}  Define ${\rm sub}_X :L^X\times L^X\longrightarrow L$ by
 $${\rm sub}_X(A,B)=\bigwedge\limits_{x\in X}A(x)\rightarrow B(x)\ (\forall A, B\in L^X).$$
 Then ${\rm sub}_X$ is an $L$-order on $L^X$, which is called the  {\em   inclusion $L$-order} on $L^X$.
If the background set is clear, then we always drop the subscript $X$ to be ${\rm sub}$.

\end{ex}

\vskip 6pt

Define ${\uparrow}x$ and ${\downarrow}x$ respectively by
${\uparrow}x(y)=e(x,y)$, ${\downarrow}x(y)=e(y,x)\ (\forall x,y\in P)$. An $L$-subset $S\in L^P$ is
called a {\it lower set} (resp., an {\it upper set}) if $S(x)\otimes e(y,x)\leq S(y)$ (resp., $S(x)\otimes e(x,y)\leq S(y)$) for all $x,y\in P$.
Clearly,  ${\downarrow}x$ (resp.,  ${\uparrow}x$) is a lower (resp., an upper) set for every  $x\in P$.
\begin{dn} {\em(\cite{PartI})} Let $P$ be an $L$-ordered set. An element $x_0$ is called a {\it supremum} of an $L$-subset $A$ of $P$, in symbols $x_0=\sqcup A$, if

{\rm (1)} $\forall x\in P,\ A(x)\leq e(x,x_0)$;

{\rm (2)} $\forall y\in P,\ \bigwedge\limits_{x\in P}A(x)\rightarrow e(x,
y)\leq e(x_0,y)$.
\\An element $x_0$ is called a {\it infimum} of an $L$-subset $A$ of $P$, in symbols $x_0=\sqcap A$,
if

{\rm (3)}  $\forall x\in P,\ A(x)\leq e(x_0,x)$;

{\rm (4)} $\forall y\in P,\ \bigwedge\limits_{x\in P}A(x)\rightarrow
e(y,x)\leq e(y,x_0)$.
\end{dn}
It is easy to verify that  if the supremum (resp., infimum) of an $L$-ordered exists, it must be unique. The following provides an equivalent definition of  supremum and infimum, respectively.

\begin{lm}{\em(\cite{PartI})} Let $P$ be an $L$-ordered set. An element $x_0\in P$ is the  supremum (resp., infimum) of $A\in L^P$, i.e., $x_0=\sqcup A$ (resp., $x_0=\sqcap A$), iff
\vskip 6pt
\centerline{$e(x_0,y)={\rm sub}(A,{\downarrow}y)\ (\forall y\in P)\ (resp., e(y,x_0)={\rm sub}(A,{\uparrow}y)\ (\forall y\in P).$}
\vskip 3pt
\end{lm}
For ease of reference, we present some simple yet useful lemmas. They belong to the folklore in the theory of $L$-orders.
\begin{lm}\label{lm-sup-leq} Let $P$ be an $L$-ordered set and let $A, B\in L^P$. If $A\leq B$, then $e(\sqcup A, \sqcup B)\geq u$.
\end{lm}

\begin{lm}
For every $L$-subset $\mathcal{A}$ in the $L$-ordered set $(L^X, {\rm sub})$, the supremum (resp., infimum) of $\mathcal{A}$ exists;
that is
\vskip3pt
\centerline{$\sqcup \mathcal{A}=\bigvee_{A\in L^X}\mathcal{A}(A)\otimes A$ (resp., $\sqcap\mathcal{A}=\bigwedge_{A\in L^X}\mathcal{A}(A)\to A$).}

\end{lm}

\begin{lm}\label{lm-join-uni} Let $\tau\subseteq L^X$ and $\mathcal{B}\in L^{\tau}$. If $\bigvee_{A\in \tau}\mathcal{B}(A)\otimes A\in \tau$, then in $(\tau, {\rm sub})$, $\bigvee_{A\in \tau}\mathcal{B}(A)\otimes A$ is the supremum of $\mathcal{B}$; that is, $\sqcup \mathcal{B}=\bigvee_{A\in \tau}\mathcal{B}(A)\otimes A$.
\end{lm}

\begin{dn}\label{dn-dir-set} {\em(\cite{PartI})} Let $P$ be an $L$-ordered set. An $L$-subset $D$ of $P$ is said to be  {\em directed}  if

{\em (D1)} $\bigvee\limits_{x\in P}D(x)\geq u$;

{\em (D2)} $\forall x,y\in P,\ D(x)\otimes D(y)\leq\bigvee\limits_{z\in
P}D(z)\otimes e(x,z)\otimes e(y,z)$.\\
 A  directed $L$-subset $I\in L^P$ is called
an ideal if it is also a lower set. The set of all
directed $L$-subsets $($resp.,  ideals$)$ of $P$ is denoted by
$\mathcal{D}_L(P)$ $($resp., $\mathcal{I}_L(P))$.
\end{dn}
% Let $P$ be an $L$-ordered set. A nonempty $L$-subset $D\in L^P$ is said to be {\it directed} if
%$D(x)\otimes D(y)\leq\bigvee\limits_{z\in P}D(z)\otimes e(x,z)\otimes e(y,z)$ $(\forall x,y\in P)$.
%An $L$-subset $I\in L^P$ is called an {\it ideal} of $P$ if it is a directed lower set.
%Denote by $\mathcal{D}_L(P)\ ($resp., $Idl_L(P))$ the set of all directed $L$-subsets (resp., ideals) of $P$.
\begin{dn}{\em(\cite{PartI})}
An $L$-ordered set $P$ is called an $L$-{\rm dcpo} if every directed $L$-subset has a supremum,
or equivalently, every ideal has a supremum.
\end{dn}

\begin{dn}{\em(\cite{PartI})}
A map $f:P\longrightarrow Q$ between two $L$-ordered sets is said to be {\rm Scott continuous} if for every $D\in\mathcal{D}_L(P)$,
 $f(\sqcup D)=\sqcup f^\rightarrow(D)$; or equivalently, for every $I\in \mathcal{I}_L(P)$,  $f(\sqcup I)=\sqcup f^\rightarrow(I)$.
\end{dn}

\subsection{fuzzy domains}
\begin{dn}\label{dn-con-po}{\em(\cite{PartI})}
Let $P$ be an  $L$-dcpo. For all $x\in P$, define $\dda x\in L^P$ by
\begin{align*}
\dda x(y)&=\bigwedge\limits_{I\in\mathcal{I}_L(P)}e(x,\sqcup I)\rightarrow I(y)\\
&=\bigwedge\limits_{D\in \mathcal{D}_L(P)}e(x,\sqcup D)\rightarrow (\bigvee_{d\in P}D(d)\otimes e(y,d)).
\end{align*}
 A  $L$-dcpo $P$ is called a  {\em continuous $L$-dcpo}, or an {\em $L$-domain},
if $\dda x$ is directed and $\sqcup{\dda}x=x$ for every $x\in P$.
\end{dn}

 \begin{lm}{\rm (\cite{PartI})}\label{lm-aux-way} Let $P$ be an $L$-dcpo. Then
$(1)$ $\forall x\in P$, $\dda x\leq {\downarrow} x$;
$(2)$ $\forall x, u, v, y\in P$, $e(u, x)\otimes \dda y(x)\otimes e(y, v)\leq\dda v(u)$.
Thus for every $x\in P$, $\dda x$ is a lower set.

\end{lm}

 \begin{lm} {\rm (\cite{PartI})}\label{lm-con-int} If $P$ is a continuous $L$-dcpo, then $\dda y(x)=\bigvee_{z\in P}\dda y(z)\otimes\dda z(x)$ for all $x,y\in P$.
\end{lm}

The following lemma plays a key role in the representation theorem for continuous  $L$-dcpos.
\begin{lm}\label{lm-use-rep} Let $P$ be a continuous $L$-dcpo. Then $(\{\dda x\mid x\in P\}, {\rm sub})\cong(P, e)$.
\end{lm}
\noindent{\bf Proof.}
Define $f: P\longrightarrow\Psi(P)$ by $f(x)=\dda x$ for all $x\in P$. Since $P$ is a continuous $L$-dcpo,  $f$ is a bijection.
For all $x, y\in P$,   since $\sqcup\dda x=x$ and $\dda y\leq {\downarrow} y$, we have
$$e(x,y)={\rm sub} (\dda x, {\downarrow} y)\geq {\rm sub} (\dda x, \dda y)={\rm sub} (f(x), f(y)).$$
By Lemma \ref{lm-aux-way}, we have $e(x,y)\otimes\dda x(t)\leq \dda y(t)$ for all $t\in P$. Therefore,
$$e(x, y)\leq  {\rm sub} (\dda x, \dda y)={\rm sub} (f(x), f(y)).$$
Thus $(\{\dda x\mid x\in P\}, {\rm sub})\cong(P, e)$. \hfill$\Box$
\begin{dn}{\rm (\cite{yao-zhao-entc})}\label{dn-alg-dcp}
 Let $P$ be an  $L$-dcpo.
 The element $x\in P$ is said to be  {\rm compact} if  $\dda x(x)\geq u$. The set of all compact elements of $P$ is denoted by $K(P)$.
 Define $k(x)\in L^{P}$ by $k(x)(y)=e(y, x)$ for $y\in K(P)$ and
otherwise $0$.
%$$
%\ \ k(x)(y)=\left\{\begin{array}{ll}e(y, x),& y\in K(P);\\
%\ \ 0,&  y\notin K(P),
%\end{array}\right.
%$$
If for every $x\in P$, $k(x)\in \mathcal{D}_{L}(P)$ and $x=\sqcup k(x)$, then $P$ is called an {\em algebraic $L$-dcpo}.
\end{dn}

\begin{rk}\label{rk-yao-alg}{\rm (\cite{yao-zhao-entc})}
{\rm (1)} $\dda x(x)\geq u$ iff $I(x)=e(x, \sqcup I)$ $(\forall I\in \mathcal{I}_L(P))$ iff

$$e(x, \sqcup D)=\bigvee_{d\in P}D(d)\otimes e(x, d)\ (\forall D\in \mathcal{D}_{L}(P)).$$

{\rm (2)} An  $L$-dcpo $P$ is algebraic if and only if for every $x\in P$,
 $k(x)|_{K(P)}\in \mathcal{D}_L(K(P))$ and $\sqcup k(x)=x$, where $k(x)|_{K(P)}(y)=e(y, x)$ for every $y\in K(P)$.

%{\rm (3)} Clearly,  for every algebraic $L$-dcpo $P$, $k(x)\leq \dda x$. It follows directly from Lemma \ref{lm-base} that every  algebraic $L$-dcpo $P$ is a continuous $L$-dcpo. But in the following, we will show that an algebraic $L$-dcpo may not be a continuous $L$-dcpo.

\end{rk}
Similar to Lemma \ref{lm-use-rep}, we have the following useful lemma.
\begin{lm}\label{lm2-use-rep} Let $P$ be an  algebraic $L$-dcpo. Then
$(\{k(x)\mid x\in P\}, {\rm sub})\cong(P, e)$.
\end{lm}

\subsection{Categorical equivalent}
For  a category $\mathcal{C}$, it is customary to denote the class of objects of $\mathcal{C}$
by $ob(\mathcal{C})$ and  the class of morphisms of $\mathcal{C}$ by $Mor(\mathcal{C})$. For all $ A, B\in ob(\mathcal{C})$, we use the symbol $Mor_{\mathcal{C}}(A, B)$ to denote the class of morphisms from $A$ to $B$.

\begin{lm}{\em(\cite{Well})}\label{lm-cat-eqv} Let $\mathcal{C}, \mathcal{D}$ be two categories. Then $\mathcal{C}$ and $\mathcal{D}$ are  equivalent
if and only if there exists a functor $F: \mathcal{C}\longrightarrow \mathcal{D}$ satisfies the following conditions:

$(1)$ $F$ is  full, namely, for all $A, B\in ob(\mathcal{C})$, $g\in Mor_{\mathcal{D}}(F(A), F(B))$,  there exists  $f\in Mor_{\mathcal{C}}(A, B)$ such that
$F(f)=g$;

$(2)$ $F$ is faithful, namely, for all $A, B\in ob(\mathcal{C})$, $f, g\in Mor_{\mathcal{C}}(A, B)$, if $f\neq g$, then $F(f)\neq F(g)$;

$(3)$ for all $B\in ob(\mathcal{D})$, there exists $A\in ob(\mathcal{C})$ such that $F(A)\cong B$.
\end{lm}

\section{The representation of continuous $L$-dcpos}

In \cite{Wang-Li-studia}, Wang and Li  introduced the notion of generalized closure spaces and used it to
represent classical  domains. The definition of a generalized closure
space merely drops expansivity. In order to extend the representation theory for  domains to  the quantale-valued setting, the first step is
to extend the notion of  generalized closure spaces to the quantale-valued setting.

\begin{dn} A {\em generalized $L$-closure space} $(X, \langle\cdot\rangle)$ is a set equipped with
an operation $\langle\cdot\rangle : L^{X}\longrightarrow L^{X}$ satisfying

{\em (GC1)} ${\rm sub}(A, B)\leq {\rm sub} (\langle A\rangle, \langle B\rangle)$ for all $A, B\in L^X$;

{\em (GC2)} $\langle\langle A\rangle\rangle\leq\langle A\rangle$ for all $A\in L^X$.

\end{dn}
As usual, we often write $X$, instead of $(X, \langle\cdot\rangle)$ for a generalized $L$-closure space. To avoid confusion, we sometimes use $\langle\cdot\rangle_X$ to emphasize that  the operator is defined  on $X$. If the background set is clear,  we always drop the subscript.
\begin{pn}\label{pn-gecl-ope} Let $(X, \langle\cdot\rangle)$ be a generalized $L$-closure space. Then

{\em (1)} ${\rm sub} (A, \langle B\rangle)\leq {\rm sub}  (\langle A\rangle, \langle B\rangle)$ for all $A, B\in L^X$.

{\em (2)} $\langle A\rangle(x)\leq {\rm sub} (\langle u_x\rangle, \langle A\rangle)$ for all $A\in L^X$ and $x\in X$.

\end{pn}
\noindent{\bf Proof.}
(1) By Conditions (GC1) and (GC2), we have
$${\rm sub} (A, \langle B\rangle)\leq {\rm sub} (\langle A\rangle,    \langle\langle B\rangle\rangle)\leq{\rm sub} (\langle A\rangle,    \langle B\rangle). $$
(2) By Lemma \ref{lm-resi-lat}(2) and Part (1),  $\langle A\rangle(x)= {\rm sub} (u_x,\langle A\rangle)\leq {\rm sub} (\langle u_x\rangle, \langle A\rangle)$.

\hfill$\Box$

\begin{dn}\label{dn-int} A generalized $L$-closure space $(X, \langle\cdot\rangle)$ is said to be  {\em interpolative} if  for all $x\in X$, $\langle\cdot\rangle$  satisfies the following
conditions:

{\em (IT1)} $\bigvee_{t\in X}\langle u_x \rangle(t)\geq u$;

{\em (IT2)} $\langle u_x \rangle(y)\leq\bigvee_{t\in X}\langle u_x \rangle(t)\otimes \langle u_t \rangle(y)$ for all $y\in X$;

{\em (IT3)} $\langle u_x\rangle(a)\otimes \langle u_x\rangle(b)\leq\bigvee_{t\in X}\langle u_x \rangle(t)\otimes \langle u_t \rangle(a)\otimes \langle u_t \rangle(b)$
for all $ a, b\in X$.

\end{dn}

\begin{dn}\label{dn-dir-clo} Let $(X, \langle\cdot\rangle)$ be an interpolative generalized $L$-closure space. An $L$-subset $U$ of $X$ is
called a {\em directed closed set} of $(X, \langle\cdot\rangle)$ if $U$ satisfies the following conditions:

{\em (DC1)} $\bigvee_{x\in X}U(x)\geq u$;

{\em (DC2)}  $U(x)\leq {\rm sub} (\langle u_x\rangle, U)$ for all $x\in X$;

{\em (DC3)}  $U(x)\leq \bigvee_{y\in X}U(y)\otimes \langle u_y\rangle(x)$ for all $x\in X$;

{\em (DC4)} $U(x)\otimes U(y)\leq \bigvee_{z\in X}U(z)\otimes \langle u_z\rangle(x)\otimes  \langle u_z\rangle(y)$ for all $x, y\in X$.

 In this paper, denote the collection of all the   directed closed sets of $(X, \langle\cdot\rangle)$ by $\mathfrak{C}(X)$.
\end{dn}

\begin{pn}\label{pn-dir-clo} Let $(X, \langle\cdot\rangle)$ be an interpolative generalized $L$-closure space.

$(1)$ $\forall x\in X$, $\langle u_x\rangle\in \mathfrak{C}(X)$;

%$(2)$ $\forall E\in \mathfrak{C}(X)$ and $A\in L^X$, $sub(A, E)\leq sub (\langle A\rangle, E)$;

$(2)$  For each directed $L$-subset $\mathcal{D}$ of $(\mathfrak{C}(X), {\rm sub})$, $\bigvee_{U\in \mathfrak{C}(X)}\mathcal{D}(U)\otimes U\in \mathfrak{C}(X)$.

\end{pn}
\noindent{\bf Proof.}
(1) It is routine to check by Proposition \ref{pn-gecl-ope}(2) and Conditions  (IT1)-(IT3) in Definition \ref{dn-int}.

(2) We  show that $\bigvee_{U\in \mathfrak{C}(X)}\mathcal{D}(U)\otimes U\in  \mathfrak{C}(X)$
satisfies
 Conditions (DC1)-(DC4) in Definition \ref{dn-dir-clo}.

  (DC1): By $\mathcal{D}$ satisfying Condition (D1) and $U$ satisfying Condition (DC1), we have
 $$\bigvee_{x\in X}\bigvee_{U\in \mathfrak{C}(X)} \mathcal{D}(U)\otimes U(x)=\bigvee_{U\in \mathfrak{C}(X)} \mathcal{D}(U)\otimes \bigvee_{x\in X}U(x)\geq u\otimes u=u.                   $$
(DC2): For all $x, y\in X$, by Condition (DC2), it holds that
\begin{align*}
\ &(\bigvee_{U\in \mathfrak{C}(X)} \mathcal{D}(U))\otimes U(x))\otimes\langle  u_x\rangle(y)\\
\leq\ &\bigvee_{U\in \mathfrak{C}(X)} \mathcal{D}(U)\otimes {\rm sub}(\langle  u_x\rangle, U)\otimes\langle  u_x\rangle(y)\\
\leq\ & \bigvee_{U\in \mathfrak{C}(X)} \mathcal{D}(U)\otimes U(y).
\end{align*}
By the arbitrariness of $y$,  $\bigvee_{U\in \mathfrak{C}(X)}\mathcal{D}(U)\otimes U(x)\leq {\rm sub}(\langle  u_x\rangle, \bigvee_{U\in \mathfrak{C}(X)}\mathcal{D}(U)\otimes U)$.

(DC3): For all $x\in X$, by Condition (DC3), it holds that
\begin{align*}
\ & \bigvee_{U\in \mathfrak{C}(X)} \mathcal{D}(U)\otimes U(x)\\
\leq \ & \bigvee_{U\in \mathfrak{C}(X)} \mathcal{D}(U)\otimes \bigvee_{y\in X}U(y)\otimes \langle  u_y\rangle(x)\\
\leq \ &\bigvee_{y\in X}(\bigvee_{U\in \mathfrak{C}(X)} \mathcal{D}(U)\otimes U(y))\otimes \langle  u_y\rangle(x).
\end{align*}
This shows that $\bigvee_{U\in \mathfrak{C}(X)} \mathcal{D}(U)\otimes U$ satisfies Condition (DC3).

(DC4): For all $a, b\in X$, by Condition (DC4) and the directedness of $\mathcal{D}$, we have
 \begin{align*}
\ &(\bigvee_{U\in \mathfrak{C}(X)} \mathcal{D}(U)\otimes U(a))\otimes (\bigvee_{U\in \mathfrak{C}(X)} \mathcal{D}(U)\otimes U(b))\\
=\ &\bigvee_{U_1, U_2\in \mathfrak{C}(X)} \mathcal{D}(U_1)\otimes\mathcal{D}(U_2)\otimes U_1(a)\otimes U_2(b)\\
\leq\ & \bigvee_{U_1, U_2,U_3\in \mathfrak{C}(X)}\mathcal{D}(U_3)\otimes {\rm sub} (U_1, U_3)\otimes {\rm sub} (U_2, U_3)\otimes U_1(a)\otimes U_2(b)\\
\leq\ &  \bigvee_{U_3\in \mathfrak{C}(X)}\mathcal{D}(U_3)\otimes U_3(a)\otimes U_3(b)\\
\leq\ &  \bigvee_{U_3\in \mathfrak{C}(X)}\mathcal{D}(U_3)\otimes \bigvee_{c\in X}U_3(c)\otimes \langle  u_c\rangle(a)\otimes \langle  u_c\rangle(b)\\
=\ &\bigvee_{c\in X}(\bigvee_{U_3\in \mathfrak{C}(X)}\mathcal{D}(U_3)\otimes U_3(c))\otimes \langle  u_c\rangle(a)\otimes \langle  u_c\rangle(b).
\end{align*}
Thus,  $\bigvee_{U\in \mathfrak{C}(X)}\mathcal{D}(U)\otimes U$ satisfies Condition (DC4).

To sum up, $\bigvee_{U\in \mathfrak{C}(X)}\mathcal{D}(U)\otimes U\in \mathfrak{C}(X)$.
\hfill$\Box$

By Lemma \ref{lm-join-uni} and Proposition \ref{pn-dir-clo}(2), we know  that for each interpolative generalized $L$-closure space $(X, \langle\cdot\rangle)$,
 $(\mathfrak{C}(X), {\rm sub})$ is an $L$-dcpo, and the  supremum of every directed $L$-subset $\mathcal{D}$ in   $(\mathfrak{C}(X), {\rm sub})$ is
  exactly $\bigvee_{U\in \mathfrak{C}(X)} \mathcal{D}(U)\otimes U$; that is $\sqcup\mathcal{D}= \bigvee_{U\in \mathfrak{C}(X)} \mathcal{D}(U)\otimes U$.

 In this paper, we write $\Psi(X)=\{\langle  u_x\rangle\mid x\in X\}$. Below, we provide some characterizations of directed closed sets.
\begin{pn}\label{pn-char-dir} Let $(X, \langle\cdot\rangle)$ be an interpolative generalized $L$-closure space. Then the following statements are equivalent.

$(1)$ $U\in \mathfrak{C}(X)$;

$(2)$ $U=\bigvee_{x\in X}U(x)\otimes \langle  u_x\rangle$;

$(3)$ There exists a directed $L$-subset $\mathcal{D}$ of $\Psi(X)$ such that
\vskip 8pt
\centerline{$\bigvee_{x\in X }\mathcal{D}(\langle u_x\rangle)\otimes \langle u_x\rangle=U$.}
\end{pn}

\noindent{\bf Proof.}
$(1)\Rightarrow (2)$: It  is obvious by  Conditions (DC2), (DC3) in Definition \ref{dn-dir-clo}.

$(2)\Rightarrow (3)$: Define $\mathcal{D}: \Psi(X)\longrightarrow L$ by
\vskip 6pt
\centerline{$\mathcal{D}(\langle  u_x\rangle)=\bigvee_{\langle  u_t\rangle=\langle  u_x\rangle}U(t)$}
\vskip 6pt
\noindent for all $x\in X$. We only need to show that  $\mathcal{D}$ is directed. By Condition (DC1),   $\mathcal{D}$ is nonempty.
 For all $x, y\in X$, by Proposition \ref{pn-gecl-ope}(2), we have
\begin{align*}
\mathcal{D}(\langle  u_x\rangle)\otimes \mathcal{D}(\langle  u_y\rangle)&=(\bigvee_{\langle  u_t\rangle=\langle  u_x\rangle}U(t))\otimes
(\bigvee_{\langle  u_w\rangle=\langle  u_y\rangle}U(w))\\
&=\bigvee \{U(t)\otimes U(w)\mid \langle  u_t\rangle=\langle  u_x\rangle,\langle  u_w\rangle=\langle  u_y\rangle\}\\
&\leq\bigvee \{\bigvee_{z\in X}U(z)\otimes \langle u_z\rangle(t)\otimes \langle u_z\rangle(w)\mid \langle  u_t\rangle=\langle  u_x\rangle,\langle  u_w\rangle=\langle u_y\rangle\}\\
&\leq \bigvee_{z\in X}U(z)\otimes {\rm sub} (\langle u_x\rangle,\langle u_z\rangle)\otimes {\rm sub} (\langle u_y\rangle,\langle u_z\rangle)\\
&\leq \bigvee_{z\in X}\mathcal{D}(\langle  u_z\rangle)\otimes {\rm sub} (\langle u_x\rangle,\langle u_z\rangle)\otimes {\rm sub} (\langle u_y\rangle,\langle u_z\rangle).
\end{align*}
Thus $\mathcal{D}$ is directed.

$(3)\Rightarrow (1)$: Define  $\mathcal{D}^{*}: \mathfrak{C}(X)\longrightarrow L$  by $\mathcal{D}^{*}(U)=\mathcal{D}(U)$ for $U\in \Psi(X)$;  $\mathcal{D}^{*}(U)=0$ for $U\notin \Psi(X)$.
% $$
%\ \ \mathcal{D}^{*}(U)=\left\{\begin{array}{ll}\mathcal{D}(U),&  U\in \Psi(X);\\
%\  \ 0,& else.
%\end{array}\right.
%$$
It is clear that $\mathcal{D}^{*}$ is a directed $L$-subset of $\mathfrak{C}(X)$.  By Proposition \ref{pn-dir-clo}(2),
 $\sqcup\mathcal{D}^{*}=\bigvee_{x\in X }\mathcal{D}(\langle u_x\rangle)\otimes \langle u_x\rangle=U$;
thus $U\in \mathfrak{C}(X)$.
\hfill$\Box$
%$(1)\Rightarrow (4)$: It  is obvious by the condition (2), (3) in Definition ?.

%$(4)\Rightarrow (1)$

We  generalize Su and Li's result (see \cite[Proposition 4.3]{LiQG1}) to the quantale-valued setting.
\begin{lm}\label{lm-base} Let $P$ be an $L$-dcpo and $x\in P$. If there exists a directed $L$-subset $D\in L^P$ such that $D\leq \dda x$ and $\sqcup D=x$, then $\dda x$ is  directed and $x=\sqcup\dda x$.
\end{lm}
\noindent{\bf Proof.} Since $D$ is directed and $D\leq \dda x$, we have
$\bigvee_{a\in P}\dda x(a)\geq D(a)\geq u.$
Thus $\dda x$ is nonempty.
For all $a, b\in P$, by $\sqcup D=x$, it holds that
\begin{align*}
\ &\dda x(a)\otimes \dda x(b)\\
\leq \ & (e(x, \sqcup D)\rightarrow\bigvee_{d_1\in X}D(d_1)\otimes e(a, d_1))\otimes e(x, \sqcup D)\rightarrow(\bigvee_{d_2\in X}D(d_2)\otimes e(b, d_2))\\
\leq \ & (u\rightarrow\bigvee_{d_1\in X}D(d_1)\otimes e(a, d_1))\otimes (u\rightarrow\bigvee_{d_2\in X}D(d_2)\otimes e(b, d_2))\\
=\ & \bigvee_{d_1\in X}D(d_1)\otimes e(a, d_1)\otimes\bigvee_{d_2\in X}D(d_2)\otimes e(b, d_2)\\
=\ & \bigvee_{d_1, d_2\in X}D(d_1)\otimes D(d_2)\otimes e(a, d_1)\otimes e(b, d_2)\\
\leq\ &\bigvee_{d_1, d_2, d_3\in X}D(d_3)\otimes e(d_1, d_3)\otimes e(d_2, d_3)\otimes e(a, d_1)\otimes e(b, d_2)\\
\leq\ &\bigvee_{ d_3\in X}D(d_3)\otimes e(a, d_3)\otimes e(b, d_3)\\
\leq\ &\bigvee_{ d_3\in X}\dda x(d_3)\otimes e(a, d_3)\otimes e(b, d_3).
\end{align*}
Thus, $\dda x$  is directed.

By $D\leq \dda x\leq {\downarrow} x$ and Lemma \ref{lm-sup-leq}, we have $e(\sqcup D, \sqcup\dda x)\geq u$ and $e(\sqcup \dda x, \sqcup{\downarrow} x)\geq u$. Since $\sqcup D=x$ and  $\sqcup{\downarrow} x=x$, we have $\sqcup\dda x=x$.
\hfill$\Box$

\begin{tm}\label{tm-gcl-dom} The $L$-ordered set  $(\mathfrak{C}(X), {\rm sub} )$ is a continuous $L$-dcpo.
\end{tm}
\noindent{\bf Proof.} By Proposition \ref{pn-dir-clo}(2),  we know that $(\mathfrak{C}(X), {\rm sub} )$ is an $L$-dcpo. Let $U\in \mathfrak{C}(X)$.
Define $\mathcal{D}_{U}: \mathfrak{C}(X)\longrightarrow L$ by
 $$\mathcal{D}_{U}(V)=\bigvee_{\langle u_x\rangle =V}U(x)$$
 for all $V\in \mathfrak{C}(X)$.
It follows from Proposition \ref{pn-char-dir} that  $\mathcal{D}_{U}$ is directed $L$-subset of  $(\mathfrak{C}(X), {\rm sub}) $ and $\sqcup\mathcal{D}_{U}=U$.

 We now prove that $\mathcal{D}_U\leq \dda U$. For every $x\in X$ and every directed $L$-subset $\mathcal{D}$ of $\mathfrak{C}(X)$, by (DC2), it holds that
\begin{align*}
U(x)\otimes {\rm sub} (U, \sqcup\mathcal{D})\leq \sqcup\mathcal{D}(x)&=\bigvee_{C\in \mathfrak{C}(X)}\mathcal{D}(C)\otimes C(x)\\
&\leq \bigvee_{C\in \mathfrak{C}(X)}\mathcal{D}(C)\otimes {\rm sub} (\langle u_x\rangle, C).
\end{align*}
Thus, for every $x\in X$, we have
\begin{align*}
\mathcal{D}_{U}(\langle u_x\rangle)&=\bigvee_{\langle u_x\rangle=\langle u_t\rangle}U(t)\\
&\leq\bigwedge_{\mathcal{D}\in \mathcal{D}_{L}(\mathfrak{C}(X))}{\rm sub} (U, \sqcup\mathcal{D})\longrightarrow
\bigvee_{C\in \mathfrak{C}(X)}\mathcal{D}(C)\otimes {\rm sub} (\langle u_x\rangle, C).
\end{align*}
This shows that $\mathcal{D}_U\leq \dda U$. By Lemma \ref{lm-base}  and  the arbitrariness of $U$,
we know that $(\mathfrak{C}(X), {\rm sub} )$ is a continuous $L$-dcpo. \hfill$\Box$
%
% $$
%\ \ \mathcal{D}_{E}(U)=\left\{\begin{array}{ll}\bigvee_{\langle u_x\rangle =U}E(x),&  U\in \Psi(X);\\
%\  \ 0,& else.
%\end{array}\right.
%$$

Theorem \ref{tm-gcl-dom} shows that every interpolative generalized $L$-closure space
can induce a continuous  $L$-dcpo through its directed closed sets.
Conversely, the following example shows that every continuous $L$-dcpo $(P, e)$ can
induce an  interpolative generalized $L$-closure space.
%Conversely, given a continuous $L$-dcpo $(P, e)$, define an operator
%$$\langle\cdot\rangle: L^P\longrightarrow L^P$$
% by
% $$\langle A\rangle=\bigvee_{x\in X}A(x)\otimes\dda x$$
%for all $A\in L^X$.
\begin{ex}\label{ex-dom-gcl} Let $P$ be a continuous $L$-dcpo.  Define an operator
$$\langle\cdot\rangle: L^P\longrightarrow L^P$$
 by
 $$\langle A\rangle=\bigvee_{x\in X}A(x)\otimes\dda x$$
for all $A\in L^X$. Then $(P, \langle\cdot\rangle)$ is an  interpolative generalized $L$-closure space.
\end{ex}
\noindent{\bf Proof.} We first prove that  $(P, \langle\cdot\rangle)$ is a generalized $L$-closure space.

 (GC1):
For all $A, B\in L^P$ and for each $t\in P$, it holds that
\begin{align*}
{\rm sub} (A, B)\otimes \langle A\rangle(t)&=\bigwedge_{x\in X}A(x)\rightarrow B(x))\otimes\bigvee_{x\in X}A(x)\otimes\dda x(t)\\
&\leq \bigvee_{x\in X}(A(x)\rightarrow B(x))\otimes A(x)\otimes\dda x(t)\\
&\leq \bigvee_{x\in X}B(x)\otimes\dda x(t)=\langle B\rangle(t).
\end{align*}
Thus, ${\rm sub} (A, B)\leq {\rm sub} (\langle A\rangle, \langle B\rangle)$.

(GC2): For all $A\in L^P$ and for each $t\in P$, by Lemma \ref{lm-aux-way}, it holds that
\begin{align*}
\langle\langle A\rangle\rangle(t)&=\bigvee_{x\in X}\langle A\rangle(x)\otimes\dda x(t)\\
&=\bigvee_{x, y\in X} A(y)\otimes\dda y(x)\otimes\dda x(t)\\
&\leq \bigvee_{ y\in X} A(y)\otimes\dda y(t)=\langle A\rangle(t).
\end {align*}
Thus $\langle\langle A\rangle\rangle\leq \langle A\rangle$.
This shows that  $(P, \langle\cdot\rangle)$ is a generalized $L$-closure space. We now prove that  $(P, \langle\cdot\rangle)$ is  interpolative. For all   $x\in P$, we have the following.

(IT1): Since $P$ is a continuous $L$-dcpo,   we have $\bigvee_{t\in P}\langle  u_x\rangle(t)=\bigvee_{t\in P}\dda x (t)\geq u$.

(IT2): It follows from Lemma \ref{lm-con-int} that
$$\langle  u_x\rangle(t)=\dda x(t)=\bigvee_{y\in X}\dda x(y)\otimes \dda y(t)=\bigvee_{y\in X}\langle  u_x\rangle(y)\otimes \langle  u_y\rangle(t)$$
for all $t\in P$.

(IT3): For $a, b\in P$,  by Lemmas \ref{lm-aux-way} and \ref{lm-con-int},  it holds that
\begin{align*}
\langle  u_x\rangle(a)\otimes \langle  u_x\rangle(b)&=\dda x(a)\otimes \dda x (b)\\
&=\bigvee_{c\in P}\dda x(c)\otimes \dda c(a)\otimes\bigvee_{d\in P} \dda x (d)\otimes \dda d(b)\\
&=\bigvee_{c, d\in P}\dda x(c)\otimes \dda x (d)\otimes \dda c(a)\otimes \dda d(b)\\
&\leq \bigvee_{c, d\in P}(\bigvee_{k\in P}\dda x(k)\otimes e(c, k)\otimes e(d, k))\otimes \dda c(a)\otimes \dda d(b)\\
&\leq \bigvee_{k\in P}\dda x(k)\otimes \dda k(a)\otimes \dda k(b)\\
&=\bigvee_{k\in P}\langle  u_x\rangle(k)\otimes \langle  u_k\rangle(a)\otimes \langle  u_k\rangle(b).
\end{align*}
To sum up, $(P, \langle\cdot\rangle)$ is an  interpolative generalized $L$-closure space. \hfill $\Box$

\begin{pn}\label{pn-dom-clo}Let $P$ be a continuous $L$-dcpo.  Then $\mathfrak{C}(P)=\{\dda x\mid x\in P\}$.
\end{pn}
\noindent{\bf Proof.}  By proposition \ref{pn-dir-clo}(1), $\{\dda x\mid x\in P\}\subseteq\mathfrak{C}(P)$ is clear. Conversely,  let $I\in \mathfrak{C}(P)$,
by Proposition \ref{pn-char-dir}, there exists a directed $L$-subset $\mathcal{D}$ of $\Psi(P)$ such that $I=\bigvee_{x\in X}\mathcal{D}(\dda x)\otimes \dda x$.
Define $D:P\longrightarrow L$ by $D(x)=\mathcal{D}(\dda x)$. Therefore, $I=\bigvee_{x\in X}D(x)\otimes \dda x$.  By Lemma \ref{lm-use-rep}, we know that  $D$ is a directed $L$-subset of $P$. Thus $\sqcup D$ exists. We claim that $I=\dda \sqcup D$. In fact, by Lemmas \ref {lm-aux-way} and \ref{lm-con-int}, we have
\begin{align*}
\dda \sqcup D(x)&=\bigvee_{y\in P}\dda \sqcup D(y)\otimes \dda y(x)\\
&=\bigvee_{d\in P}D(d)\otimes e(y,d)\otimes \dda y(x)\\
&= \bigvee_{d\in P}D(d)\otimes  \dda d(x)\\
&=I(x).
\end{align*}
This shows that $\mathfrak{C}(P)\subseteq\{\dda x\mid x\in P\}$. Thus $\mathfrak{C}(P)=\{\dda x\mid x\in P\}$.
\hfill $\Box$

%$$I=\bigvee_{x\in P}\mathcal{D}(\langle u_x\rangle)\otimes \langle u_x\rangle=\mathcal{D}(\dda x)\otimes \dda x.$$
We now arrive at the main conclusion of this section.
\begin{tm}\label{tm-rep-con} {\em (Representation theorem I: for continuous $L$-dcpos)} An $L$-ordered set $(P, e)$ is a continuous $L$-dcpo if and only if
there exists  an interpolative generalized $L$-closure space $(X, \langle\cdot\rangle)$ such that $(\mathfrak{C}(X), {\rm sub})\cong (P, e)$.
\end{tm}
\noindent{\bf Proof.}
 {\bf Sufficiency.} It follows directly from Theorem  \ref{tm-gcl-dom}.

{\bf Necessity.}
It  follows from Lemma  \ref{lm-use-rep}, Example \ref{ex-dom-gcl} and Proposition  \ref{pn-dom-clo}.
\hfill $\Box$
\section{The representation of algebraic $L$-dcpos}
In this section, we explore  representations of algebraic $L$-dcpos. Let us first review the research history of representing algebraic domains in the classical setting.
A well-known fact is the Stone-type duality between algebraic closure spaces/convex spaces and algebraic lattices, demonstrating that algebraic lattices can be represented by algebraic closure spaces/convex spaces (see \cite{Erne-stone,shen-zhao-shi}).
To overcome the limitation that algebraic closure spaces cannot directly represent algebraic domains, Guo and Li introduced the notion of F-augmented closure spaces, which are triples, to successfully provide representations for algebraic domains \cite{Guo-Li}.
In \cite{Wu-Guo-Li}, Wu, Guo, and Li utilized a specific subfamily of closed sets in algebraic closure spaces to provide representations for algebraic domains.
 In \cite{Wang-Li-studia}, Wang and Li utilized a special type of interpolative generalized closure spaces to represent algebraic domains.
 Recently, in \cite{Wu-Xu-FIE}, Wu and Xu simplified the conditions for F-augmented closure spaces and algebraic closure spaces mentioned in \cite{Guo-Li,Wu-Guo-Li}, extending the representation of algebraic domains to the framework of classical closure space.
It is worth noting that there are no essential differences in the methods presented in \cite{Guo-Li,Wu-Guo-Li,Wu-Xu-FIE}.
Another aim of this section is to
 discuss the relationship between the representation methods  for algebraic domains in \cite{Wang-Li-studia}
 and those in \cite{Guo-Li,Wu-Guo-Li,Wu-Xu-FIE} within the framework of quantale-valued settings.

\begin{dn}\label{dn-l-clo}{\em (\cite{Bel-2001})} An  $L$-closure space $(X, \langle\cdot\rangle)$ is a set equipped with
an operation $\langle\cdot\rangle : L^{X}\longrightarrow L^{X}$ that satisfies the following conditions:

{\em (LC1)} $A\leq \langle A\rangle$ for all $A\in L^X$;

{\em (LC2)}  ${\rm sub} (A, B)\leq {\rm sub} (\langle A\rangle, \langle B\rangle)$ for all $A, B\in L^X$;

{\em (LC3)}  $\langle\langle A\rangle\rangle=\langle A\rangle$ for all $A\in L^X$.
\end{dn}

\begin{pn} Let $(X, \langle\cdot\rangle)$ be a generalized  $L$-closure space. Then $X$ is an $L$-closure space if and only if for every $x\in X$, $\langle u_x\rangle(x)\geq u$.
\end{pn}
\noindent{\bf Proof.} The part of necessity is clear. We only need to show the part of sufficiency.
For ever $A\in L^X$, by (GC1), it holds that
$$ A(x)={\rm sub} (u_x, A)\leq {\rm sub} (\langle u_x\rangle,  \langle A\rangle)\leq  \langle A\rangle(x).$$
This shows that $A\leq \langle A\rangle$.\hfill $\Box$

It is evident that  by adding the Condition (LC1)  to the definition of a generalized  $L$-closure space, we can obtain the definition of an $L$-closure space.
 Notice that Condition (LC1) ensures that  $L$-closure space automatically satisfies  Conditions (IT1)-(IT3). Thus, each $L$-closure space is an interpolative generalized
 $L$-closure space.

To show that  each $L$-closure space can induce an algebraic $L$-dcpo, we first present a useful lemma.
\begin{lm}\label{lm-bas-alg}Let $P$ be an  $L$-dcpo and $x\in P$. If there exists a directed L-subset $D$ such that $D\leq k(x)$ (cf. Definition \ref{dn-alg-dcp}) and $\sqcup D=x$, then $ k(x)$ is directed and $x=\sqcup k(x)$.
\end{lm}
\noindent{\bf Proof.} By $D\leq k(x)\leq {\downarrow} x$,  it is clear that   $\sqcup k(x)=x$ and $\bigvee_{y\in P} k(x)(y)\geq u$.
For every $y_{1}, y_2\in K(P)$,  we have
\begin{align*}
k(x)(y_1)\otimes k(x)(y_2)
&=e(y_1, x)\otimes e(y_2, x)\\
&=e(y_1, \sqcup D)\otimes e(y_2, \sqcup D)\\
&= (\bigvee_{d\in P}D(d)\otimes e(y_1, d))\otimes (\bigvee_{d\in P}D(d)\otimes e(y_2, d))\\
&=\bigvee_{d_1, d_2\in P}(D(d_1)\otimes D(d_2)\otimes e(y_1, d_1)\otimes e(y_2, d_2))\\
&=\bigvee_{d_1, d_2\in P}\bigvee_{d\in P}(D(d)\otimes e(d_1, d)\otimes e(d_2, d)\otimes e(y_1, d_1)\otimes e(y_2, d_2))\\
&\leq\bigvee_{d\in P} k(x)(d)\otimes e(y_1, d)\otimes e(y_2, d).
\end{align*}
Thus $ k(x)$ is  directed.
\hfill$\Box$
\begin{tm} \label{tm-clo-ald} Let $(X, \langle\cdot\rangle)$  be an $L$-closure space. Then $(\mathfrak{C}(X), {\rm sub} )$ is an algebraic $L$-dcpo.
\end{tm}
\noindent{\bf Proof.} By Theorem \ref{tm-gcl-dom}, $(\mathfrak{C}(X), {\rm sub} )$ is a continuous  $L$-dcpo. Let $U\in \mathfrak{C}(X)$. For every $x\in X$ and for every directed $L$-subset $\mathcal{D}$ of $\mathfrak{C}(X)$, by Conditions (LC1) and (DC2), it holds that
\begin{align*}
{\rm sub} (\langle u_x\rangle, \sqcup \mathcal{D})&\leq \sqcup \mathcal{D}(x)=\bigvee_{U\in \mathfrak{C}(X)}\mathcal{D}(U)\otimes U(x)\\
&\leq \bigvee_{U\in \mathfrak{C}(X)}\mathcal{D}(U)\otimes {\rm sub} (\langle u_x\rangle, U).
\end{align*}
This shows that $\langle u_x\rangle\in K(\mathfrak{C}(X))$. Thus $\Psi(X)=\{\langle  u_x\rangle\mid x\in X\}\subseteq K(\mathfrak{C}(X))$.
Let  $\mathcal{D}_U$ be the directed $L$-subset of  $\mathfrak{C}(X)$ defined in Theorem \ref{tm-gcl-dom}; that is
$\mathcal{D}_U(V)=\bigvee_{V=\langle u_x\rangle}U(x)$
for all $V\in \mathfrak{C}(X)$. Notice that $\Psi(X)\subseteq K(\mathfrak{C}(X))$
and for every $x\in X$,
$$\mathcal{D}_U(\langle u_x\rangle)=\bigvee_{\langle u_t\rangle=\langle u_x\rangle}U(t)\leq {\rm sub}(\langle u_x\rangle, U).$$
We have  $\mathcal{D}_U\leq k(U)$. Thus by Lemma \ref{lm-bas-alg},
$(\mathfrak{C}(X), {\rm sub} )$ is an algebraic $L$-dcpo.

\hfill$\Box$

Theorem \ref{tm-clo-ald} shows that every  $L$-closure space
can induce an  algebraic  $L$-dcpo by making use of its directed closed sets.
Conversely, the following example shows that every algebraic $L$-dcpo $(P, e)$ can
induce an $L$-closure space in a natural way.

\begin{ex}\label{ex-algdom-cl}
Let $(P, e)$ be an algebraic $L$-dcpo. Define $\langle\cdot\rangle: L^{K(P)}\longrightarrow L^{K(P)} $ by
$$\langle A\rangle= {\downarrow}A:=\bigvee_{x\in K(P)}A(x)\otimes {\downarrow} x,$$
for all $A\in L^{K(P)}$. It is claer that $(K(P),\langle\cdot\rangle)$
is an $L$-closure space.
\end{ex}
\begin{pn}\label{pn-clo-diclo}Let $(P, e)$ be an algebraic $L$-dcpo. Then
$\mathfrak{C}(K(P))=\{{\downarrow} x|_{K(P)}\mid x\in P\}$.
\end{pn}
\noindent{\bf Proof.} By Definition \ref{dn-alg-dcp}, it is routine to check that  $\{{\downarrow} x|_{K(P)}\mid x\in P\}\subseteq \mathfrak{C}(K(P))$.
Conversely, Let $U\in \mathfrak{C}(K(P))$. Then by Proposition \ref{pn-char-dir}(3),
there exists a directed $L$-subset $\mathcal{D}$ of ${\Psi(K(P))}$ such that
$$U=\bigvee_{k\in K(P)}\mathcal{D}(\langle u_k\rangle)\otimes\langle u_k\rangle=\bigvee_{k\in K(P)}\mathcal{D}({\downarrow} k|_{K(P)})\otimes{\downarrow} k|_{K(P)}.$$
Define $D:P\longrightarrow L$ by
 $$
\ \ D(x)=\left\{\begin{array}{ll}\mathcal{D}(\langle u_x\rangle),&  x\in K(P);\\
\  \ 0,& else.
\end{array}\right.
$$
 It is clear that $D$ is nonempty. Notice that for each $x\in K(P)$, $\langle u_x\rangle={\downarrow} x|_{K(P)}(=k(x)|_{K(P)})$. By Lemma \ref{lm2-use-rep}
 and the fact that $\mathcal{D}$ is directed, we know that  $D$ is also a directed $L$-subset of $P$. Thus $\sqcup D$ exists. We claim that $U={\downarrow} \sqcup D|_{K(P)}$. In fact,
for each $x\in K(P)$,
\begin{align*}
{\downarrow} \sqcup D|_{K(P)}(x)&=e(x,  \sqcup D)=\bigvee_{d\in P}D(d)\otimes e(x, d)\\
&=\bigvee_{d\in K(P)}D(d)\otimes e(x, d)\\
&=\bigvee_{d\in K(P)}\mathcal{D}({\downarrow} d|_{K(P)})\otimes {\downarrow} d|_{K(P)}(x)=U(x).
\end{align*}
This shows that $(\mathfrak{C}(K(P))\subseteq\{{\downarrow} x|_{K(P)}\mid x\in P\}$, as desired.
\hfill$\Box$

%\begin{align*}
%D(x_1)\otimes D(x_2)&=\mathcal{D}(\langle u_{x_1}\rangle)\otimes\mathcal{D}(\langle u_{x_2}\rangle)\\
%&\leq\bigvee_{x_3\in K(P)}\mathcal{D}(\langle u_{x_3}\rangle)\otimes {\rm sub}(\langle u_{x_1}\rangle, \langle u_{x_3}\rangle)\otimes {\rm sub}(\langle u_{x_2}\rangle, \langle u_{x_3}\rangle)\\
%&=\bigvee_{x_3\in K(P)}D(x_3)\otimes e(x_1, x_3)\otimes e(x_2, x_3)
%\end{align*}
We now obtain the first main result of this section, which is a quantale-valued type of the representation theorem for algebraic domains as seen in  \cite{Guo-Li,Wu-Guo-Li,Wu-Xu-FIE}.
\begin{tm}\label{tm-rep1-alg} {\em (Representation theorem II: for algebraic $L$-dcpos)} An $L$-ordered set $(P, e)$ is an algebraic $L$-dcpo if and only if there exists
an  $L$-closure space $(X, \langle\cdot\rangle)$ such that $(\mathfrak{C}(X), {\rm sub} )\cong (P, e)$.
\end{tm}
\noindent{\bf Proof.} {\bf Sufficiency.} It follows directly from Theorem  \ref{tm-clo-ald}.

{\bf Necessity.} Let $(P, e)$ be an algebraic $L$-dcpo. Notice that
$(\{{\downarrow} x|_{K(P)}\mid x\in P\}, {\rm sub})\cong(\{ k(x)\mid x\in P\}, {\rm sub})$.
 By Example \ref{ex-algdom-cl}  and  Lemma \ref{lm2-use-rep}, we have $(\mathfrak{C}(K(P)), {\rm sub} )\cong (P, e)$.
\hfill $\Box$

In the remainder of this section, we assume that  $L$ is a {\bf commutative integral quantale}, i.e.,  $u=1$. We will investigate an alternative representation for algebraic $L$-dcpos.

\begin{dn} Let $(X, \langle\cdot\rangle)$ be an interpolative generalized $L$-closure space and $Y\subseteq X$. If for every $x\in X$, it holds that

%{\rm(DS1)} $\bigvee_{y\in X}\langle1_x\rangle(y)=1$;

{\rm(DS)} $\langle1_x\rangle(a)\leq \bigvee_{y\in Y}\langle1_x\rangle(y)\otimes\langle1_y\rangle(a)$ for all $a\in X$,

%{\rm(DS3)} $\langle1_x\rangle(a)\otimes \langle1_x\rangle(b)\leq \bigvee_{y\in Y}\langle1_x\rangle(y)\otimes\langle1_y\rangle(a)\otimes\langle1_y\rangle(b)$ for all $a, b\in X$,

then it is routine to check that  $(Y, \langle\cdot\rangle|_Y)$ is an interpolative generalized $L$-closure space, called a {\em dense subspace of $(X, \langle\cdot\rangle)$}, where $\langle\cdot\rangle|_Y$ is the restriction of $\langle\cdot\rangle$ on $Y$.
\end{dn}

\begin{pn}\label{pn-clos-iso} Let $(Y, \langle\cdot\rangle|_Y)$  be a dense subspace of an  interpolative generalized $L$-closure space $(X, \langle\cdot\rangle)$.
Then $(\mathfrak{C}(Y), {\rm sub} )\cong(\mathfrak{C}(X), {\rm sub} )$.
\end{pn}
\noindent{\bf Proof.} We divide this proof into three steps.

{\bf Step 1.} For each $E\in \mathfrak{C}(X)$, we prove that  $E|_Y\in \mathfrak{C}(Y)$.

(DC1):
\begin{align*}
\bigvee_{x\in X}E(x)&=\bigvee_{x, y\in X}E(y)\otimes \langle1_y\rangle(x)&(\text{ Proposition} \ref{pn-char-dir})\\
&\leq\bigvee_{x, y\in X}E(y)\otimes (\bigvee_{z\in Y}\langle1_y\rangle(z)\otimes \langle1_z\rangle(x))&(\text{DS})\\
&=\bigvee_{x, y\in X, z\in Y}E(y)\otimes\langle1_y\rangle(z)\otimes \langle1_z\rangle(x)\\
&\leq \bigvee_{x\in X, z\in Y}E(z)\otimes \langle1_z\rangle(x)\\
&= \bigvee_{z\in Y}E(z)\otimes \bigvee_{x\in X}\langle1_z\rangle(x)\\
&= \bigvee_{z\in Y}E(z).&(\text{IT1})
\end{align*}
By $\bigvee_{x\in X}E(x)=1$, we have $\bigvee_{z\in Y}E|_Y(z)=\bigvee_{z\in Y}E(z)=1$.

(DC2):  For every $y\in Y$, by $E\in \mathfrak{C}(X)$, it is clear that  $E|_Y(y)\leq{\rm sub}(\langle1_y\rangle|_Y, E|_Y)$.

(DC3): For every $x\in Y$,  by $E\in \mathfrak{C}(X)$, we have
\begin{align*}
 E|_Y(x)&=E(x)\leq \bigvee_{y\in X}E(y)\otimes \langle1_y\rangle(x)\\
 &\leq\bigvee_{y\in X}E(y)\otimes \bigvee_{z\in Y}\langle1_y\rangle(z)\otimes\langle1_z\rangle(x)&(\text{DS})\\
&\leq \bigvee_{z\in Y}E(z)\otimes \langle1_z\rangle(x)\\
&=\bigvee_{z\in Y}E|_Y(z)\otimes \langle1_z\rangle|_Y(x)
\end{align*}
Thus $E|_Y(x)\leq E|_Y(z)\otimes \langle1_z\rangle|_Y(x)$.

(DC4): For every $a, b\in Y$,
\begin{align*}
 E|_Y(a)\otimes E|_Y(b)&=E(a)\otimes E(b)\\
 &\leq\bigvee_{c\in X}E(c)\otimes \langle1_c\rangle(a)\otimes\langle1_c\rangle(b)&(\text{DC4})\\
&\leq \bigvee_{c\in X, d\in X}E(c)\otimes \langle1_c\rangle(d)\otimes\langle1_d\rangle(a)\otimes\langle1_d\rangle(b)&(\text{IT3})\\
&\leq \bigvee_{c\in X, d\in X}E(c)\otimes (\bigvee_{y\in Y}\langle1_c\rangle(y)\otimes\langle1_y\rangle(d)) \otimes\langle1_d\rangle(a)\otimes\langle1_d\rangle(b)&(\text{DS})\\
&\leq \bigvee_{ y\in Y}E(y)\otimes \langle1_y\rangle(a)\otimes\langle1_y\rangle(b)\\
&=\bigvee_{ y\in Y}E|_Y(y)\otimes \langle1_y\rangle|_Y(a)\otimes\langle1_y\rangle_Y(b).
\end{align*}

{\bf Step 2.} We prove that ${\rm sub} _Y({E_1}|_Y, {E_2}|_Y)={\rm sub} _X({E_1}, {E_2})$.

We only need to show ${\rm sub} _Y({E_1}|_Y, {E_2}|_Y)\leq {\rm sub} _X({E_1}, {E_2})$. For all $x\in X$,
\begin{align*}
{\rm sub} _Y({E_1}|_Y, {E_2}|_Y)\otimes E_1(x)&\leq (\bigwedge_{y\in Y}E_1(y)\to E_2(y))\otimes \bigvee_{z\in X}E_1(z)\otimes\langle1_z\rangle(x)\\
&\leq (\bigwedge_{y\in Y}E_1(y)\to E_2(y))\otimes \bigvee_{z\in X, t\in Y}E_1(z)\otimes\langle1_z\rangle(t)\otimes \langle1_t\rangle(x)\\
&\leq \bigvee_{t\in Y}(\bigwedge_{y\in Y}E_1(y)\to E_2(y))\otimes E_1(t)\otimes \langle1_t\rangle(x)\\
&\leq \bigvee_{t\in Y}(E_1(t)\to E_2(t))\otimes E_1(t)\otimes \langle1_t\rangle(x)\\
&\leq  E_2(t)\otimes \langle1_t\rangle(x)
\leq E_2(x).
\end{align*}
By the arbitrariness of $x$, we have ${\rm sub} _Y({E_1}|_Y, {E_2}|_Y)\leq {\rm sub} _X({E_1}, {E_2})$.

{\bf Step 3.} We prove that $\mathfrak{C}(Y)=\{E|_Y\mid E\in \mathfrak{C}(X)\}$.

By Step 1, we only need to prove that $\mathfrak{C}(Y)\subseteq\{E|_Y\mid E\in \mathfrak{C}(X)\}$. Let $S\in \mathfrak{C}(Y)$. By Proposition \ref{pn-char-dir}(2),
$S=\bigvee_{y\in Y}S(y)\otimes \langle1_y\rangle|_Y$.
Define $\mathcal{D}:\Psi(X)\longrightarrow L$ by
$$\mathcal{D}(\langle1_x\rangle)=\bigvee\{S(y)\mid \langle1_x\rangle=\langle1_y\rangle, y\in Y\}.$$
 It is routine to check that $\mathcal{D}$ is a directed $L$-subset of $\Psi(X)$. Define  $\mathcal{D}^{*}: \mathfrak{C}(X)\longrightarrow L$  by $\mathcal{D}^{*}(U)=\mathcal{D}(U)$ for $U\in \Psi(X)$;  $\mathcal{D}^{*}(U)=0$ for $U\notin \Psi(X)$. It is easy to see that $\mathcal{D}^{*}$
  is also a directed $L$-subset of $\mathfrak{C}(X)$ and $\sqcup\mathcal{D}^{*}=\bigvee_{y\in Y}S(y)\otimes \langle1_y\rangle$. By Proposition \ref{pn-dir-clo}(2), $\bigvee_{y\in Y}S(y)\otimes \langle1_y\rangle\in \mathfrak{C}(X)$.
It follows that $$S=(\bigvee_{y\in Y}S(y)\otimes \langle1_y\rangle)|_{Y}\in \{E|_Y\mid E\in \mathfrak{C}(X)\}.$$
Based on the above steps, we have $(\mathfrak{C}(Y), {\rm sub} )\cong(\mathfrak{C}(X), {\rm sub} )$.\hfill $\Box$

The verification of (DC1) in Step 1 of the above proof   depends on $L$ being integral. This is why we need $L$ is integral.
%Define $\mathcal{D}:\Psi(Y)\longrightarrow L$ by
%$$\mathcal{D}(\langle1_y\rangle|_Y)=\bigvee\{S(t)\mid \langle1_y\rangle=\langle1_t\rangle, t\in Y\}$$
%for all $y\in Y$.
%$L$-subset of $\Psi(Y)$.
%Notice that
%$$\sqcup i^{\rightarrow}(\mathcal{D})=\bigvee_{y\in Y}(\bigvee_{\langle1_t\rangle=\langle1_y\rangle, t\in Y}S(t))\otimes \langle1_y\rangle=\bigvee_{y\in Y}S(y)\otimes \langle1_y\rangle.$$
\begin{ex}\label{ex-dens-sub} Let $(P, e)$ be an algebraic $L$-dcpo. Then $(K(P), \langle\cdot\rangle|_{K(P)})$ is a dense subspace of $(P, \langle\cdot\rangle)$ (cf. Example \ref{ex-dom-gcl}),
 where $\langle\cdot\rangle|_{K(P)}$ refers to the restriction of   $ \langle\cdot\rangle$ on  $K(P)$. It is worth noting that
 $(K(P), \langle\cdot\rangle|_{K(P)})$ is exactly the $L$-closure space defined  in Example \ref{ex-algdom-cl}.
\end{ex}
\noindent{\bf Proof.}  For each  $x, a\in P$,
\begin{align*}
\langle1_x\rangle(a)=\dda x(a)=\dda \sqcup k(x)(a)&\leq \bigvee_{y\in P}k(x)(y)\otimes e(a,y)\\
&=\bigvee_{y\in K(P)}\dda x(y)\otimes e(a,y)\\
&=\bigvee_{y\in K(P)}\langle1_x\rangle(y)\otimes\langle1_y\rangle(a).
\end{align*}

%(%DS3): For all $a, b\in P$, by the proof of the Part (DS2), it holds that
%\begin{align*}
%\langle1_x\rangle(a)\otimes \langle1_x\rangle(b)&=\dda x(a)\otimes \dda x(b)\\
%&\leq\bigvee_{y_1\in K(P)}\dda x(y_1)\otimes e(a, y_1)\otimes\bigvee_{y_2\in K(P)}\dda x(y_2)\otimes e(b, y_2) \\
%&=\bigvee_{y_1,y_2\in K(P)}\dda x(y_1)\otimes \dda x(y_2)\otimes e(a, y_1)\otimes e(b, y_2)\\
%&\leq \bigvee_{y_1,y_2\in K(P)}\bigvee_{y_3\in K(P)}\dda x(y_3)\otimes e(y_1, y_3)\otimes e(y_2, y_3)\otimes e(a, y_1)\otimes e(b, y_2)\\
%&\leq \bigvee_{y\in K(P)}\dda x(y)\otimes e(a, y)\otimes e(b, y)\\
%&=\bigvee_{y\in K(P)}\langle1_x\rangle (y)\otimes \langle1_y\rangle(a)\otimes  \langle1_y\rangle(b)
%\end{align*}
Thus, $(K(P), \langle\cdot\rangle|_{K(P)})$ is a dense subspace of $(P, \langle\cdot\rangle)$
\hfill $\Box$

We now obtain an alternative representation for algebraic $L$-dcpos, marking the second main conclusion of this section.
It is worth noting again that the truth value table
$L$  here is a commutative  integral quantale.

\begin{tm}\label{wang-alg-fuzzy}{\em (Representation Theorem III: for algebraic $L$-dcpos)}
An  $L$-ordered set $(P, e)$ is an algebraic $L$-dcpo  if and only if there exists
an interpolative generalized $L$-closure space $(X, \langle\cdot\rangle)$  with an  $L$-closure  space as
a dense subspace such that $(\mathfrak{C}(X), {\rm sub} )\cong (P, e)$.
\end{tm}
\noindent{\bf Proof.}
{\bf Sufficiency.} Suppose the $L$-closure space $Y$ is a dense subspace of $X$.
By Proposition \ref{pn-clos-iso}, we have $(\mathfrak{C}(X), {\rm sub} )\cong (\mathfrak{C}(Y), {\rm sub} )$. Therefore $(\mathfrak{C}(Y), {\rm sub} )\cong (P, e)$.
By Theorem \ref{tm-clo-ald}, $(P, e)$ is an algebraic $L$-dcpo.

{\bf Necessity.}  Since every algebraic $L$-dcpo is  a continuous $L$-dcpo, by Theorem \ref{tm-rep-con},
 we have $(\mathfrak{C}(P),{\rm sub} )\cong(P,e)$.
By Example \ref{ex-dens-sub}, the $L$-closure space $(K(P), \langle\cdot\rangle|_{K(P)})$ is a dense subspace of $(P, \langle\cdot\rangle)$.
\hfill $\Box$

Theorem \ref{wang-alg-fuzzy} provides a quantale-valued version  of Wang and Li's representation  for algebraic domains (see \cite[Theorem 3.20]{Wang-Li-semi}).
However, we achieve this representation  from the perspective of a more general dense subspace. Proposition \ref{pn-clos-iso}
 reveals the essential role of dense subspaces, thus   establishing a link between  the methods of \cite{Guo-Li,Wu-Guo-Li,Wu-Xu-FIE}
 and those of \cite{Wang-Li-studia} within the framework of the quantale-valued setting.

\section{Approximable $L$-relations}

In classical domain theory, Scott continuous mappings between two posets are the most typical and often serve as morphisms in categories of domains.
In order to represent Scott continuous mappings
between domains, Scott introduced the notion of an approximable
mapping between his information systems \cite{Scott-82}.  Unlike usual mappings between two sets, Scott's approximable mappings for information
systems are binary relations that satisfy  certain axioms.

Recently,  Scott-type  approximable mappings for closure spaces, along with many other similar binary relations, have been widely considered as a convenient approach to reconstruct Scott continuous mappings (see \cite{Guo-Li, Li-Wang-Yao, Wu-Guo-Li, Wang-Li-studia, Wu-Xu-FIE, Yao1-Li, Yao2-Li}).
In this section, we aim to explore the relationship between continuous $L$-dcpos and   interpolative generalized $L$-closure spaces from a categorical perspective.
 We first extend the notion of  approximable mappings for interpolative generalized closure spaces to that  for   quantale-valued ones.

%Following Scott's work, Larsen and Winskel
%made use of  Scott¡¯s approximable mappings  to represent Scott
%continuous mappings between Scott domains.
% In [24], Scott introduced the notion of an approximable
%mapping between his information systems. Soon after, Larsen and Winskel
%proved that Scott¡¯s approximable mappings can be used to represent Scott
%continuous functions between Scott domains [18]. Unlike Ern¡äe¡¯s continuous
%functions for closure spaces, Scott¡¯s approximable mappings for information
%systems are relations that satisfy some axioms. Scott¡¯s approximable mappings, as same as many similar relations, have been widely considered to be
%a good approach to capturing Scott continuous functions [1,12,13,25,32].
%Following the idea of Scott¡¯s approximable mappings for information systems,
%we now introduce the notion of an approximable mapping for IG-closure spaces. We will see that our approximable mappings between
%IG-closure spaces can induce a category equivalent to that of domains with
%Scott continuous functions.

\begin{dn}\label{dn-app-rel} Let $(X, \langle\cdot\rangle_X)$, $(Y, \langle\cdot\rangle_Y)$ be two interpolative
generalized $L$-closure spaces and let $\Theta\in L^{X\times Y}$ be an $L$-relation.
$\Theta$ is called an {\em approximable $L$-relation} from $(X, \langle\cdot\rangle_X)$ to $(Y, \langle\cdot\rangle_Y)$ if it satisfies the following conditions:

{\em (AP1)} $\bigvee_{y\in Y}\Theta(x, y)\geq u$ for all $x\in X$;

{\em (AP2)} $\langle u_{x^{\prime}}\rangle(x)\otimes \Theta(x, y)\leq \Theta(x^{\prime}, y)$ for all $x, x^{\prime}\in X$ and $y\in Y$;

{\em (AP3)} $\Theta(x, y)\otimes \langle u_{y}\rangle( y^{\prime})\leq\Theta(x, y^{\prime}) $ for all $x\in X$ and $y, y^{\prime}\in Y$;

{\em (AP4)}  $\Theta(x, y)\leq \bigvee_{x^{\prime}\in X, y^{\prime}\in Y}\langle u_{x}\rangle(x^{\prime})\otimes \Theta(x^{\prime}, y^{\prime})
\otimes\langle u_{y^{\prime}}\rangle( y)$ for all $x\in X$ and $ y\in Y$;

{\em (AP5)} $\Theta(x, y_1)\otimes \Theta(x, y_2)\leq
\bigvee_{y_3\in Y}\Theta(x, y_3)\otimes\langle u_{y_3}\rangle(y_1)\otimes\langle u_{y_3}\rangle(y_2)$ for all $x\in X $and  $y_1, y_2\in Y$.

\end{dn}

\begin{rk}\label{rk-app-rel} {\em (1)} If $(X, \langle\cdot\rangle_X)$ and $(Y, \langle\cdot\rangle_Y)$  are two $L$-closure spaces, then Condition (AP4) holds automatically.

{\em (2)} By the conditions in Definition \ref{dn-app-rel}, we get the following useful equalities.
For all $x\in X$ and $y\in Y$,
\begin{align*}
\Theta(x, y)&=\bigvee_{x^{\prime}\in X}\langle u_{x}\rangle(x^{\prime})\otimes \Theta(x^{\prime}, y)\\
&=\bigvee_{ y^{\prime}\in Y} \Theta(x, y^{\prime})\otimes\langle u_{y^{\prime}}\rangle( y)\\
&=\bigvee_{x^{\prime}\in X, y^{\prime}\in Y}\langle u_{x}\rangle(x^{\prime})\otimes \Theta(x^{\prime}, y^{\prime})\otimes\langle u_{y^{\prime}}\rangle( y).
\end{align*}
\end{rk}

\begin{dn}\label{dn-ide-app}Let $(X, \langle\cdot\rangle)$ be an interpolative  generalized $L$-closure space. Define $id_X\in L^{X\times X}$ by
$$id_X(x, y)=\langle u_x\rangle(y)\ (\forall x, y\in X).$$
Then it is easy to check that $id_X$ is an $L$-approximable relation from $X$ to itself.
Let $\Theta$ be an $L$-approximable relation from  $X$ to  $Y$. By Remark \ref{rk-app-rel}(2), we have
 $$\Theta\circ id_X=id_X\circ\Theta=\Theta.$$
So,  $id_X$  is called {\em the identity  $L$-approximable relation} on $X$.
\end{dn}
\begin{pn}\label{pn-com-app}  Let $X$, $Y$ and $Z$ be interpolative generalized $L$-closure spaces and let  $\mathrel{\Theta}\in L^{X\times Y}$ and $\mathrel{\Upsilon}\subseteq L^{Y\times Z}$ be approximable $L$-relations. Then the composition  $\mathrel{\Upsilon}\circ\mathrel{\Theta}$  is an approximable $L$-relation from $X$ to $Z$.
\end{pn}
\noindent{\bf Proof.}
It is clear  that  $\mathrel{\Upsilon}\circ\mathrel{\Theta}$
satisfies Condition  (AP1) in Definition \ref{dn-app-rel}.

(AP2): For all $x, x^{\prime}\in X$ and for all $z\in Z$, it holds that
\begin{align*}
&\langle u_{x^{\prime}}\rangle (x)\otimes {\Upsilon}\circ\Theta(x, z)\\
=&\langle u_{x^{\prime}}\rangle(x)\otimes \bigvee_{y\in Y}{\Theta}(x, y)\otimes \Upsilon(y,z)\\
=&\bigvee_{y\in Y}{\Theta}(x^{\prime}, y)\circ {\Upsilon}(y,z)\\
=&{\Upsilon}\circ{\Theta}(x^{\prime}, z).
\end{align*}
%Thus $\langle u_{x^{\prime}}\rangle (x)\otimes {\Upsilon}\circ\Theta(x, z)={\Upsilon}\circ{\Theta}(x^{\prime}, z)$.

(AP3):  Notice that $\Upsilon$ satisfies the Condition (AP3). Similar to the verification of (AP2),
 it is routine to check that  $\mathrel{\Upsilon}\circ\mathrel{\Theta}$ satisfies Condition (AP3).

(AP4): For all  $x\in X$  and for all $ z\in Z$, by Remark \ref{rk-app-rel}(2), it holds that
\begin{align*}
\Upsilon\circ \Theta(x, z)&=\bigvee_{y\in Y}\Theta(x, y)\otimes \Upsilon(y,z)\\
&= \bigvee_{y\in Y}\bigvee_{x^{\prime}\in X}\langle u_{x}\rangle(x^{\prime})\otimes \Theta(x^{\prime}, y)\otimes \bigvee_{z^{\prime}\in Z}\Upsilon(y,z^{\prime})\otimes \langle u_{z^{\prime}}\rangle(z)\\
&=\bigvee_{y\in Y,x^{\prime}\in X, z^{\prime}\in Z}\langle u_{x}\rangle(x^{\prime})\otimes \Theta(x^{\prime}, y)\otimes \Upsilon(y,z^{\prime})\otimes \langle u_{z^{\prime}}\rangle(z)\\
&=\bigvee_{x^{\prime}\in X, z^{\prime}\in Z}\langle u_{x}\rangle(x^{\prime})\otimes \Upsilon\circ\Theta(x^{\prime}, z^{\prime})\otimes \langle u_{z^{\prime}}\rangle(z).
\end{align*}
%Thus $\Upsilon\circ \Theta(x, z)=
%\bigvee_{x^{\prime}\in X, z^{\prime}\in Z}\langle u_{x}\rangle(x^{\prime})\otimes \Upsilon\circ\Theta(x^{\prime}, z^{\prime})\otimes \langle u_{z^{\prime}}\rangle(z)$.

(AP5): For all  $x\in X$ and $z_1, z_2\in Z$, it holds that
\begin{align*}
&\Upsilon\circ \Theta(x, z_1)\otimes\Upsilon\circ \Theta(x, z_2)\\
=&\bigvee_{y_1, y_2\in Y}\Theta(x, y_1)\otimes\Theta(x, y_2)\otimes \Upsilon(y_1,z_1)\otimes\Upsilon(y_2,z_2)\\
\leq&  \bigvee_{y_1, y_2\in Y}\bigvee_{y_3\in Y}\Theta(x, y_3)\otimes \langle u_{y_3}\rangle(y_1)\otimes \langle u_{y_3}\rangle(y_2)\otimes \Upsilon(y_1,z_1)\otimes\Upsilon(y_2,z_2)\\
\leq& \bigvee_{y_3\in Y}\Theta(x, y_3)\otimes \Upsilon(y_3,z_1)\otimes\Upsilon(y_3,z_2)\\
\leq& \bigvee_{y_3\in Y, z_3\in Z}\Theta(x, y_3)\otimes \Upsilon(y_3,z_3)\otimes \langle u_{z_3}\rangle(z_1)\otimes \langle u_{z_3}\rangle(z_2)\\
=&\bigvee_{ z_3\in Z}\Upsilon\circ\Theta(x,z_3)\otimes \langle u_{z_3}\rangle(z_1)\otimes \langle u_{z_3}\rangle(z_2).
\end{align*}
%Thus $\Upsilon\circ \Theta(x, z_1)\otimes\Upsilon\circ \Theta(x, z_2)
%\leq\bigvee_{ z_3\in Z}\Upsilon\circ\Theta(x,z_3)\otimes \langle u_{z_3}\rangle(z_1)\otimes \langle u_{z_3}\rangle(z_2).$

To sum up,  $\mathrel{\Upsilon}\circ\mathrel{\Theta}$ is an approximable $L$-relation relation from $X$ to $Z$.
\hfill$\Box$

Let $\Theta$ be an  approximable $L$-relation relation from $X$ to $Y$ and $U\in \mathfrak{C}(X)$. Define $\widetilde{\Theta}(U)\in L^{Y}$ by
$$\widetilde{\Theta}(U)(y)=\bigvee_{x\in X}U(x)\otimes \Theta(x, y)\ (\forall y\in Y).$$
By Remark \ref{rk-app-rel}, it is easy to see that
$$\widetilde{\Theta}(U)(y)=\bigvee_{x\in X, y^{\prime}\in Y}U(x)\otimes \Theta(x, y^{\prime})\otimes\langle u_{y^{\prime}}\rangle(y).\ (\forall y\in Y)$$
 This fact will be useful in the following context.
\begin{pn} Let $\Theta$ be an  approximable $L$-relation relation from $(X, \langle\cdot\rangle_X)$ to $(Y, \langle\cdot\rangle_Y)$.
If $U$ is a directed closed set of $X$,
then $\widetilde{\Theta}(U)$ is a directed closed set of $Y$.
\end{pn}
\noindent{\bf Proof.}
Define $\mathcal{D}:\mathfrak{C}(Y)\longrightarrow L$ by
$\mathcal{D}(V)=\bigvee_{x\in X, \langle u_{t}\rangle=V}U(x)\otimes \Theta(x, t)$
for all $V\in \mathfrak{C}(Y)$.
Since $U$ satisfies (DC1)   and $\Theta$ satisfies (AP1), it is clear that $\mathcal{D}$  is a nonempty set of $\mathfrak{C}(Y)$. For all $y_1, y_2\in Y$,
by Conditions (DC4), (AP2) and (AP5),
\begin{align*}
&\mathcal{D}(\langle u_{y_1}\rangle)\otimes\mathcal{D}(\langle u_{y_2}\rangle)\\
=&(\bigvee_{x\in X, \langle u_{t_1}\rangle=\langle u_{y_1}\rangle}U(x)\otimes \Theta(x, t_1))\otimes (\bigvee_{x\in X, (\langle u_{t_2}\rangle=\langle u_{y_2}\rangle}U(x)\otimes \Theta(x, t_2))\\
=&\bigvee\{U(x_1)\otimes U(x_2) \otimes \Theta(x_1, t_1)\otimes \Theta(x_2, t_2)\mid x_1, x_2\in X, \langle u_{t_1}\rangle=\langle u_{y_1}\rangle, \langle u_{t_2}\rangle=\langle u_{y_2}\rangle\}\\
\leq &\bigvee\{U(x_3)\otimes \langle u_{x_3}\rangle(x_1)\otimes \langle u_{x_3}\rangle(x_2) \otimes \Theta(x_1, t_1)\otimes \Theta(x_2, t_2)\mid x_1, x_2, x_3\in X, \langle u_{t_i}\rangle=\langle u_{y_i}\rangle (i=1,2)\}\\
\leq &\bigvee\{U(x_3)\otimes \Theta(x_3, t_1)\otimes \Theta(x_3, t_2)\mid  x_3\in X, \langle u_{t_i}\rangle=\langle u_{y_i}\rangle (i=1,2)\}\\
\leq &\bigvee\{\bigvee_{y_3\in Y}U(x_3)\otimes \Theta(x_3, y_3)\otimes\langle u_{y_3}\rangle(t_1)\otimes\langle u_{y_3}\rangle(t_2) \mid  x_3\in X, \langle u_{t_i}\rangle=\langle u_{y_i}\rangle (i=1,2)\}\\
\leq &\bigvee\{\bigvee_{y_3\in Y}U(x_3)\otimes \Theta(x_3, y_3)\otimes {\rm sub} (\langle u_{t_1}\rangle, \langle u_{y_3}\rangle)\otimes {\rm sub} (\langle u_{t_2}\rangle, \langle u_{y_3}\rangle)\mid  x_3\in X, \langle u_{t_i}\rangle=\langle u_{y_i}\rangle (i=1,2)\}\\
\leq &\bigvee_{y_3\in Y}\mathcal{D}(\langle u_{y_3}\rangle)\otimes {\rm sub} (\langle u_{y_1}\rangle, \langle u_{y_3}\rangle)\otimes {\rm sub} (\langle u_{y_2}\rangle, \langle u_{y_3}\rangle).
\end{align*}
This shows that $\mathcal{D}$ is a directed $L$-subset of $\mathfrak{C}(Y)$.
Notice that
\begin{align*}
\ & \bigvee_{U\in \mathfrak{C}(X)}\mathcal{D}(U)\otimes U\\
=\ & \bigvee_{y\in Y} \mathcal{D}(\langle u_{y}\rangle)\otimes \langle u_{y}\rangle\\
=\ & \bigvee_{y\in Y}(\bigvee_{x\in X, \langle u_{t}\rangle=\langle u_{y}\rangle}U(x)\otimes \Theta(x, t))\otimes \langle u_{y}\rangle\\
=\ & \bigvee_{x\in X, y\in Y}U(x)\otimes \Theta(x, y)\otimes \langle u_{y}\rangle\\
=\ &\widetilde{\Theta}(U).
\end{align*}
By Proposition \ref{pn-dir-clo}(2), we have $\widetilde{\Theta}(U)=\sqcup\mathcal{D}\in \mathfrak{C}(Y)$.
\hfill$\Box$

The above proposition  shows that an $L$-approximable relation $\Theta$ can give rise to an assignment from directed
 closed sets of $X$  to those of $Y$. Moreover, we will see that this assignment can induce a Scott continuous mapping between related continuous $L$-dcpos.

\begin{pn}For every   approximable $L$-relation relation $\Theta$ between two interpolative generalized $L$-closure spaces $X$ and $Y$.
Define $\psi_{\Theta}: \mathfrak{C}(X)\longrightarrow\mathfrak{C}(Y)$
by $$\psi_{\Theta}(U)=\widetilde{\Theta}(U)$$
for each $U\in \mathfrak{C}(X)$.
Then $\psi_{\Theta}$ is a Scott continuous mapping from $(\mathfrak{C}(X), {\rm sub} )$ to $(\mathfrak{C}(Y), {\rm sub} )$.
\end{pn}
\noindent{\bf Proof.}  It is clear that $\psi_{\Theta}:(\mathfrak{C}(X), {\rm sub} )\longrightarrow(\mathfrak{C}(Y), {\rm sub} )$ is $L$-order-preserving.
For every directed $L$-subset $\mathcal{D}$ of $\mathfrak{C}(X)$ and for every $y\in Y$,  By Proposition \ref{pn-dir-clo}(2), we have
$$\psi_{\Theta}(\sqcup \mathcal{D})(y)=\bigvee_{x\in X}\sqcup \mathcal{D}(x)\otimes \Theta(x, y)
=\bigvee_{x\in X, U\in\mathfrak{C}(X)}\mathcal{D}(U)\otimes U(x)\otimes \Theta (x, y).$$
It is routine to check that $\psi_{\Theta}^{\rightarrow}(\mathcal{D})$ is a directed $L$-subset of $\mathfrak{C}(Y)$. By Proposition \ref{pn-dir-clo}(2), we have
\begin{align*}
\sqcup\psi_{\Theta}^{\rightarrow}(\mathcal{D})(y)&=\bigvee_{V\in\mathfrak{C}(Y)}\psi_{\Theta}^{\rightarrow}(\mathcal{D})(V)\otimes V(y)\\
&=\bigvee_{U\in \mathfrak{C}(X)}(\bigvee_{\psi_{\Theta}(U)=V}\mathcal{D}(U))\otimes V(y)\\
&=\bigvee_{U\in \mathfrak{C}(X)}\mathcal{D}(U)\otimes \psi_{\Theta}(U)(y)\\
&=\bigvee_{x\in X, U\in \mathfrak{C}(X)}\mathcal{D}(U)\otimes U(x)\otimes \Theta(x, y).
\end{align*}
This shows that $\psi_{\Theta}(\sqcup \mathcal{D})=\sqcup\psi_{\Theta}^{\rightarrow}(\mathcal{D})$;
thus $\psi_{\Theta}$
is a Scott continuous mapping from $(\mathfrak{C}(X), {\rm sub} )$ to $(\mathfrak{C}(Y), {\rm sub} )$. \hfill$\Box$

For the converse direction, we draw the following conclusion.

%given two interpolative generalized $L$-closure spaces, every  Scott continuous mapping
%between the collections of their   directed closed sets can induce a  $L$-approximable relation
\begin{pn}
 Let $ X$ and $Y$ be two interpolative generalized $L$-closure spaces and let $\psi: \mathfrak{C}(X)\longrightarrow\mathfrak{C}(Y)$ be a Scott continuous map.
Define $\mathrel{\Theta}_{\psi}\in L^{X\times Y}$ by
$$\mathrel{\Theta}_{\psi}(x, y)=\psi(\langle u_{x}\rangle)(y)$$
 for all $(x, y)\in X\times Y$.
Then $\mathrel{\Theta}_{\psi}$ is an $L$-approximable relation from $X$ to $Y$.
\end{pn}
\noindent{\bf Proof.}
It is clear that condition (AP1) holds since $\psi(\langle u_{x}\rangle)\in \mathfrak{C}(Y)$.

(AP2): For all $x, x^{\prime}\in X$ and $y\in Y$, since $\psi$ is $L$-order-preserving, we have
\begin{align*}
\langle u_{x^{\prime}}\rangle(x)\otimes \Theta_{\psi}(x, y)&\leq {\rm sub} (\langle u_{x}\rangle, \langle u_{x^{\prime}}\rangle)
\otimes \psi(\langle u_{x}\rangle)(y)\\
&\leq \psi(\langle u_{x}^{\prime}\rangle)(y)= \Theta_{\psi}(x^{\prime}, y).
\end{align*}

(AP3): For all $x\in X$ and $y, y^{\prime}\in Y$, by Condition (DC2), we have
\begin{align*}
\Theta_{\psi}(x, y)\otimes \langle u_{y}\rangle( y^{\prime})&=\psi(\langle u_{x}\rangle)(y)\otimes \langle u_{y}\rangle( y^{\prime})\\
&\leq {\rm sub} (\langle u_{y}\rangle, \psi(\langle u_{x}\rangle))\otimes \langle u_{y}\rangle( y^{\prime})\\
&\leq  \psi(\langle u_{x}\rangle)( y^{\prime})=\Theta_{\psi}(x, y^{\prime}).
\end{align*}

(AP4):  For all $x\in X$ and $y\in Y$, by Propositions \ref{pn-dir-clo}(2), \ref{pn-char-dir}(2) and $\psi$ being Scott continuous, we have
\begin{align*}
\Theta_{\psi}(x, y)=\psi(\langle u_{x}\rangle)(y)&=\bigvee_{z\in Y}\psi(\langle u_{x}\rangle)(z)\otimes\langle u_{z}\rangle(y)\\
=&\bigvee_{z\in Y}(\bigvee_{t\in X}\langle u_{x}\rangle(t)\otimes\psi(\langle u_{t}\rangle)(z))\otimes\langle u_{z}\rangle(y)\\
&\ \ \ ( \mbox{ by }\langle u_{x}\rangle=\bigvee_{t\in X}\langle u_{x}\rangle(t)\otimes \langle u_{t}\rangle )\\
=&\bigvee_{t\in X,z\in Y}\langle u_{x}\rangle(t)\otimes\Theta_{\psi}(t, z)\otimes\langle u_{z}\rangle(y).
\end{align*}

(AP5): For $x\in X$ and $y_1, y_2\in Y$, by Condition (DC4), we have
\begin{align*}
\Theta_{\psi}(x, y_1)\otimes \Theta_{\psi}(x, y_2)&=\psi(\langle u_{x}\rangle)(y_1)\otimes \psi(\langle u_{x}\rangle)(y_2)\\
&\leq\bigvee_{y_3\in Y}\psi(\langle u_{x}\rangle)(y_3)\otimes\langle u_{y_3}\rangle(y_1)\otimes\langle u_{y_3}\rangle(y_2)\\
&\leq\bigvee_{y_3\in Y} \Theta_{\psi}(x, y_3)\otimes\langle u_{y_3}\rangle(y_1)\otimes\langle u_{y_3}\rangle(y_2).
\end{align*}
To sum up, $\mathrel{\Theta}_{\psi}$ is an $L$-approximable relation from $X$ to $Y$.
\hfill$\Box$

We now  investigate how Scott continuous mappings between continuous $L$-dcpos can be represented by the notion of  approximable $L$-relations.
\begin{pn}\label{pn-one-one}
Let $\Theta\in L^{X\times Y}$ be an approximable $L$-relation two interpolative generalized $L$-closure spaces and let  $\psi: \mathfrak{C}(X)\longrightarrow\mathfrak{C}(Y)$ be a Scott continuous mapping. Then
$\Theta_{\psi_{\Theta}}=\Theta$ and  $\psi_{{\Theta}_{\psi}}=\psi$.
\end{pn}
\noindent{\bf Proof.} For all $x\in X$ and $y\in Y$ and , it holds that
$$ \Theta_{\psi_{\Theta}}(x, y)=\psi_{\Theta}(\langle u_{x}\rangle)(y)=\bigvee_{t\in X}\langle u_{x}\rangle(t) \otimes \Theta(t, y)=\Theta(x, y).$$
This shows that $\Theta_{\psi_{\Theta}}=\Theta$.

For all $x\in X$ and for all $U\in \mathfrak{C}(X)$,
define $\mathcal{D}:\mathfrak{C}(X)\longrightarrow L$ by
$$\mathcal{D}(V)=\bigvee_{\langle u_{x}\rangle=V}U(x),$$
for all $V\in\mathfrak{C}(X)$. By Propositions \ref{pn-dir-clo}(2) and \ref{pn-char-dir}(2), we know that
$\sqcup\mathcal{D}=\bigvee_{x\in X}U(x)\otimes\langle u_{x}\rangle=U.$
Thus, for all $y\in Y$, we have the following equality:
\begin{align*}
\psi_{{\Theta}_{\psi}}(U)(y)&=\bigvee_{x\in X}U(x)\otimes  {\Theta}_{\psi}(x, y)\\
&=\bigvee_{x\in X}U(x)\otimes \psi( \langle u_{x}\rangle)(y)\\
&=\sqcup  \psi^{\rightarrow}(\mathcal{D})(y)= \psi(\sqcup\mathcal{D})(y)\\
&=\psi(U)(y).
\end{align*}
This shows that $\psi_{{\Theta}_{\psi}}=\psi$.
\hfill$\Box$

%&=\sqcup i^{\rightarrow}(\psi^{\rightarrow}(\mathcal{D}))\\
By   Definition \ref{dn-ide-app} and Proposition \ref{pn-com-app}, we know  that interpolative generalized $L$-closure spaces with $L$-approximable relation as morphisms
can form a category, denoted by  $L\text{-}\mathbf{IGCS}$. We use  $L\text{-}\mathbf{CS}$ to denote the category of $L$-closure spaces,
which is a full subcategory of $L\text{-}\mathbf{IGCS}$.
 Let $L\text{-}\mathbf{CDOM}$ be the category of continuous $L$-dcpos with Scott continuous mappings as morphisms.
 Moreover, we use $L\text{-}\mathbf{AlgDOM}$ to denote the category of algebraic $L$-dcpos, which is a full subcategory of  $L\text{-}\mathbf{CDOM}$.
 We now establish categorical equivalences between relevant categories.

\begin{tm}\label{Cat-sL-sLD}
The categories $L\text{-}\mathbf{IGCS}$ (resp., $L\text{-}\mathbf{CS}$) and  $L\text{-}\mathbf{CDOM}$ (resp., $L\text{-}\mathbf{AlgDOM}$) are equivalent.
\end{tm}
\noindent{\bf Proof.}  Define $F: L\text{-}\mathbf{IGCS}\longrightarrow L\text{-}\mathbf{CDOM}$ by
$F(X)=(\mathfrak{C}(X), {\rm sub} )$ for all $ X\in ob$($L\text{-}\mathbf{IGCS}$);
$F({\Theta})=\psi_{{\Theta}}$   for all $ \mathrel{\Theta}\in$ Mor($L\text{-}\mathbf{IGCS}$).\\
We claim that $F$ is a functor from $L\text{-}\mathbf{IGCS}$ to  $L\text{-}\mathbf{CDOM}$. For all $U\in \mathfrak{C}(X)$ and $y\in X$, by Definition \ref{dn-ide-app},
  \begin{align*}
  F(id_X)(U)(y)& = \psi_{id_X}(U)(y)\\
 &=\bigvee_{x\in X}U(x)\otimes id_X(x, y)\\
 &=\bigvee_{x\in X}U(x)\otimes \langle u_{x}\rangle(y)=U(y).
\end{align*}
This shows that $F(id_X)=id_{\mathfrak{C}(X)}$.

Let $\Theta$ be an $L$-approximable relation from $X$ to $Y$ and  let $\Upsilon$  be an $L$-approximable relation from $Y$ to $Z$.
Then for any $U\in \mathfrak{C}(X)$, we have
\begin{align*}
F(\Upsilon\circ \Theta)(U)&=\psi_{\Upsilon\circ \Theta}(U)(z)\\
&=\bigvee_{x\in X}U(x)\otimes \Upsilon\circ \Theta(x, z)\\
&=\bigvee_{x\in X, y\in Y}U(x)\otimes \Theta(x, y)\otimes\Upsilon(y,z)\\
&=\bigvee_{y\in Y}\psi_{\Theta}(U)(y)\otimes\Upsilon(y,z)\\
&=\psi_{\Upsilon}(\psi_{\Theta}(U))(z).
\end{align*}
This shows that  $F(\Upsilon\circ \Theta)=F(\Upsilon)\circ F(\Theta)$. Thus  $F$ is a functor from $L\text{-}\mathbf{IGCS}$ to  $L\text{-}\mathbf{CDOM}$.

By Theorem \ref{tm-rep-con} and Proposition \ref{pn-one-one}, it is easy to see that the functor $F$ satisfies three conditions in  Lemma \ref{lm-cat-eqv};
 thus $L\text{-}\mathbf{IGCS}$ and  $L\text{-}\mathbf{CDOM}$ are equivalent.

 Since $L\text{-}\mathbf{CS}$ and  $L\text{-}\mathbf{AlgDOM}$ are
 full subcategories of $L\text{-}\mathbf{IGCS}$ and $L\text{-}\mathbf{CDOM}$ respectively,
 it follows that $L\text{-}\mathbf{CS}$ and  $L\text{-}\mathbf{AlgDOM}$  are equivalent.
 \hfill$\Box$
\section{Conclusion}

This paper explores representations of $L$-domains using interpolative generalized $L$-closure spaces for a commutative unital quantale $L$ as the truth value table. We have established a categorical equivalence between    the category of continuous $L$-dcpos (resp., algebraic $L$-dcpos) with Scott continuous mappings  and
 that of   interpolative generalized  $L$-closure spaces (resp., $L$-closure spaces) with $L$-approximable relations.
 This  opens a way to finding
representations for $L$-domain structures via $L$-closure spaces and    strengthens the connection  between
$L$-closure spaces and  $L$-ordered structures.

We have extended  the results of Wang and Li \cite{Wang-Li-studia}
to  the quantale-valued  setting, but we didn't follow their proof methods  step by
step.
 Additionally, we  explored some new results based on their result. The main differences between our work and theirs are:
 \begin{itemize}

\item In classical settings, the definition of interpolation for a generalized closure space is typically articulated with only one condition
(see \cite[Condition (In)]{Wang-Li-studia}). However, in fuzzy settings, particularly within broader quantale-valued contexts, describing interpolation through the fuzzy counterpart of Condition (In) directly is unattainable. Consequently, we had to technically divide Condition (In) into three separate conditions without losing  its essence (see Definition \ref{dn-int}). These conditions effectively capture the interpolation characteristics of generalized $L$-closure spaces.
Similarly, direct imitation of generalized directed sets in the classical setting (see \cite[Definition 3.7]{Wang-Li-studia}) to establish a counterpart in the quantale-valued setting is not feasible. This prompted us to skillfully divide the statement for generalized directed  subsets into four conditions, thereby obtaining the notion of directed closed subsets (see Definition \ref{dn-dir-clo}).
%\item  In \cite{Zhao-Zhang}, Zhao and Zhang introduce  the notions of join $\Omega$-semilattice (see  \cite[Definition 3.8]{Zhao-Zhang}) and weak join $\Omega$-semilattice (see \cite[Definition 3.16]{Zhao-Zhang}), which are two  counterparts of join-semilattice  in fuzzy setting. In the future, we can study the relationship between the Zhao and Zhang's fuzzy join-semilattices and our $L$-join-semilattices.

%\item  %In the theory of topology, studying generalizations of sobriety is  an interesting topic.
%Zhao and Fan in \cite{zhao-fan}
% introduced a weak notion of sobriety, called bounded sobriety.  Zhang and Wang in \cite{Zhang-Wang}
% extended bounded sobriety  sobriety to the setting of $\mathcal{Q}$-cotopological spaces.
% Following these work, we can introduce the notion of bounded sobriety in the framework of convex structure;
% and then study the relationship between bounded  sober $L$-convex spaces and conditional join-semilattices.
\item
The existing methods for representing algebraic domains using closure spaces can essentially be categorized into two types.
 One is     based on classical closure spaces (see \cite{Guo-Li,Wu-Guo-Li,Wu-Xu-FIE});
  the other is based on generalized closure spaces (see  \cite{Wang-Li-studia}). It is therefore natural to wonder about the relationship between these two methods.
In Section 4, the notion of dense subspace we introduced plays an
 essential role since it possesses   a family of  directed closed sets that is order-isomorphic to that of  the large space containing it (see Proposition \ref{pn-clos-iso}). Consequently, when the truth value table $L$ is commutative integral quantale, we   established a link between  the method of \cite{Guo-Li,Wu-Guo-Li,Wu-Xu-FIE}
 and that of \cite{Wang-Li-studia} to some extend.

%
%Theorem \ref{wang-alg-fuzzy} is a lattice-type  of Wang and Li's representation  for algebraic domains ( see \cite[Theorem 3.220]{Wang-Li-semi}).
%But here we provide a representation of algebraic $L$-dcpos from the perspective of a more general dense subspace. Proposition \ref{pn-clos-iso}
%actually reveals the essential role of dense subspaces. Thus we.

\end{itemize}

 %discuss the relationship between the method of representing for algebraic domains in \cite{Wang-Li-studia}
% and the method of that in \cite{Guo-Li,Wu-Guo-Li,Wu-Xu-FIE} in the framework of fuzzy setting.

%Consequently, ordered approach  and categorical approach can be mutually used in studying objects in fuzzy convex structure theory in the future.

%We have generalized Xu and Mao's  results (see \cite{XU-mao-form}),  Mislove's results  (see \cite{Mislove}) and Lawson's results (see \cite{Lawson})
% to those in the frame-valued setting, but we didn't follow their proof methods.
%  Moreover, we  explored some new results based on their result. The main differences between our work and their work are:
\section*{Acknowledgment}

This paper is supported by National Natural Science Foundation of China
(12231007, 12371462), Jiangsu Provincial Innovative and Entrepreneurial
Talent Support Plan (JSSCRC2021521).

\end{document}